\documentclass[leqno,12pt]{article}
\usepackage{latexsym}
\usepackage{amsmath}
\usepackage{amssymb}
\usepackage{bbm}
\usepackage{color}

\usepackage{graphicx}
\usepackage{epsfig}
\usepackage{hyperref} 

 0

\newtheorem{thm}{Theorem}[section]

\newtheorem{cor}[thm]{Corollary}
\newtheorem{lem}[thm]{Lemma}

\newtheorem{rem}[thm]{Remark}

\newcommand{\ea}{\end{array}}
\newcommand{\beqohne}{\begin{eqnarray*}}
\newcommand{\eeqohne}{\end{eqnarray*}}
\newcommand{\beohne}{\begin{equation*}}
\newcommand{\eeohne}{\end{equation*}}
\newcommand{\R}{\mathbb{R}}
\newcommand{\N}{\mathbb{N}}

\newcommand{\lal}{\langle\langle}
\newcommand{\rar}{\rangle\rangle}

\newcommand{\idbold}{\boldsymbol{1}}
\newcommand{\abold}{\boldsymbol{\alpha}}
\newcommand{\rbold}{\boldsymbol{\rho}}

\newcommand{\Fbold}{\boldsymbol{\Phi}}
\newcommand{\tobold}{\boldsymbol{\tau}}
\newcommand{\kbold}{\boldsymbol{\kappa}}

\def \sur#1#2{\mathrel{\mathop{\kern 0pt#1}\limits^{#2}}}

\def \w{{\tt w}}
\def \v{{\tt v}}
\def \ben{\begin{eqnarray}}
\def \een{\end{eqnarray}}
\newcommand{\tr}{\mathrm{tr}}

\newcommand{\Sr}{\mathcal{S}}

\newcommand{\GUE}{\operatorname{GUE}}

\newcommand{\LUE}{\operatorname{LUE}}
\newcommand{\JUE}{\operatorname{JUE}}
\newcommand{\KMK}{\operatorname{KMK}}

\newcommand{\dd}{\mathrm d}

\makeatletter\@addtoreset{equation}{section}\makeatother

\begin{document}

\title{Sum rules and large deviations for spectral matrix measures in the Jacobi ensemble}
\author{{\small Fabrice Gamboa}\footnote{ Universit\'e Paul Sabatier, Institut de Math\'ematiques de Toulouse,  31062-Toulouse Cedex 9, France, 
fabrice.gamboa@math.univ-toulouse.fr}
\and{\small Jan Nagel}\footnote{Eindhoven University of Technology, Department of Technology and Computer Science, 5600 MB Eindhoven, Netherlands,
  e-mail: j.h.nagel@tue.nl}
\and{\small Alain Rouault}\footnote
{Laboratoire de Math\'ematiques de Versailles, UVSQ, CNRS, Universit\'e Paris-Saclay, 78035-Versailles Cedex France, e-mail: alain.rouault@uvsq.fr}}

\maketitle

\begin{abstract}
We continue to explore the connections between large deviations for  objects coming from random matrix theory and sum rules. 
This connection was established in \cite{magicrules} for spectral measures of classical ensembles (Gauss-Hermite, Laguerre, Jacobi) and it was extended to spectral matrix measures of the Hermite and Laguerre ensemble in \cite{GaNaRomat}. In this paper, we consider the remaining case of spectral matrix measures of the Jacobi ensemble. Our main results are a large deviation principle for such measures and a sum rule for matrix measures with reference measure the Kesten-McKay law. 
As an important intermediate step, we derive the distribution of canonical moments of the matrix Jacobi ensemble.
\end{abstract}


\section{Introduction}

A probability measure on a compact subset of $\mathbb R$ or on the unit circle may be encoded by the sequence of its moments or by the coefficients of the recursion satisfied by the corresponding orthogonal polynomials. It is however not easy to relate information on the measure, (for example on its support), with information on the recursion coefficients. 
Sum rules give a way to translate between these two languages. 
Indeed, a sum rule is an identity relating a functional of the probability measure, usually in the form of a realative entropy, and a functional of its recursion coefficients. The "{\it measure side}" of the identity gives the discrepancy between the measure and a reference measure and the "{\it coefficient side}" gives the discrepancy between the correponding series of recursion coefficients.

One of the most  classical example of such a sum rule is the Szeg{\H o}-Verblunsky theorem for measures on the unit circle $\mathbb T$, see Chapter 1 of \cite{simon2}. Here, the reference measure  is the uniform measure on $\mathbb{T}$ and the coefficient side involves a sum of functions of the Verblunsky coefficients. The most famous sum rule for measures on the real line is the Killip-Simon sum rule \cite{KS03}  (see also \cite{simon2} Section 3.5). In this case, the reference measure is the semicircle distribution. In \cite{magicrules}, we gave a probabilistic interpretation of the Killip-Simon sum rule (KS-SR) and a general  strategy to construct and prove new  sum rules. The starting point is a $N\times N$ random matrix $X_N$ chosen according to the Gaussian unitarily invariant ensemble. The random spectral measure $\mu_N$ of this random matrix is then defined through its moments, by the relation 
\begin{align*}
\int x^k d\mu_N = (X_N^k)_{1,1} .  
\end{align*}
It was shown in \cite{gamboacanonical}, that as $N$ tends to infinity, the sequence $(\mu_N)_N$ satisfies a large deviation principle (LDP). The  rate function $\mathcal{I}_c$ is a functional of the recursion coefficients. Surprinsingly, this functional is exactly 
 the coefficient side of the KS-SR. Later, in \cite{magicrules}, we  gave an alternative proof of this LDP, 
with a rate function $\mathcal{I}_m$ that is 
exactly the measure side of KS-SR. 
Since a large deviation rate function is unique, this
 implies the sum rule identity $\mathcal{I}_c=\mathcal{I}_m$. 
Working with a random matrix of one of the other two classical ensembles, the Laguerre and  Jacobi ensemble, this method leads to new sum rules. 
Here the reference measures are the Marchenko-Pastur law and the Kesten-McKay law, respectively \cite{magicrules}. 
We also refer to recent interesting 
 developments of the method explored in \cite{BSZ} and \cite{breuer2018large}.

One of the ingredient to prove 
 the LDP 
 in terms of the coefficients 
is the fact that these coefficients are independent and have explicit distributions. To be more precise, it has been  shown in \cite{dumede2002}, that in the Gaussian case the coefficients are independent random variables with normal or gamma distributions. The Laguerre case has also been considered in \cite{dumede2002}. In this last frame, 
the convenient encoding is not directly by 
the 
 recursion coefficients, but by 
decomposition of them into independent variables. In \cite{Killip1}, a further decomposition is 
 shown 
 for the Jacobi ensemble. Actually these variables are the Verblunsky coefficient of the measure lifted to the unit circle, which are sometimes also called canonical moments, see the monograph \cite{dette1997theory}. 

A natural extension of scalar measures are measures with values 
 in the space of Hermitian nonnegative definite matrices. There is a rich theory of polynomials orthogonal with respect to such a matrix measure, and we refer the interested reader to \cite{sinap1996orthogonal}, \cite{duran1996orthogonal}, \cite{duran1995orthogonal} or \cite{damanik2008analytic} and referecences therein. Surprisingly, 
 sum rule identities also hold in the matrix frame.
In \cite{damal}, 
 a matricial version of KS-SR is proved 
(see also Section 4.6 of \cite{simon2}). In \cite{GaNaRomat}, we have extended our probabilistic method to the matrix case as well, and have proved an LDP 
 for random matrix valued spectral measures. This $p\times p$ measure $\Sigma_N$ is now defined by 
its matrix moments 
\begin{align} \label{defmatrixspectral}
\int x^k d\Sigma_N(x) = (X_N^k)_{i,j=1,\dots ,p} ,\qquad  k \geq 1,
\end{align}
where $X_N$ is as before a random $N\times N$ matrix and $N\geq p$. Using the explicit construction of random matrices of the Gaussian and Laguerre ensemble, it is possible to derive the distribution of the recursion coefficients of $\Sigma_N$, which are now $p\times p$ matrices, and prove an LDP 
for them, generalizing the results of \cite{dumede2002} and \cite{gamboacanonical}. Collecting these two LDPs 
 and different representations of the rate function, we obtain the matrix sum rule both for Gaussian and Laguerre cases. 
A large deviation principle for the coefficients in the matricial Jacobi case, and consequently a new sum rule, has been open so far.    

In this paper, we complete the trio of matrix 
 measures of classical ensembles by addressing 
 the Jacobi case. We prove an LDP 
for the spectral matrix measure in Theorem \ref{thm:LDPcoefficient}, which then implies the 
 new matrix
 sum rule stated in Theorem \ref{LDPSR}. 
A crucial ingredient for the proof of Theorem \ref{thm:LDPcoefficient} 
is Theorem \ref{thm:distributioncanonical}, 
where we derive the distribution of matricial canonical moments of $\Sigma_N$.
Up to our knowledge, this result is new. Actually, we have to consider for our probabilistic approach certain Hermitian versions of the canonical moments and we show that these versions are independent and each distributed as $p\times p$ matrices of the Jacobi ensemble, thereby generalizing the results of \cite{Killip1} to the matrix case.
 An additional difficulty is that the measures we need to consider are finitely supported and then are not nontrivial. In this case, many arguments 
used in the scalar case cannot be extended directly.
 The fact that there is still a one-to-one correspondence between the spectral measure $\Sigma_N$ and its canonical moments might therefore be of independent interest.

Let us explain the main obstacle 
 that so far impeded a large deviation analysis of the coefficients in the 
matrix Jacobi case. For 
 the Gaussian or Laguerre ensemble, the distribution of recursion coefficients can be derived through repeated Householder reflections applied to the full matrix $X_N$. In the Jacobi case, it seems impossible to control the effect of these tranformation on the different subblocks of $X_N$. Instead, looking at the scalar case, there are two potential strategies. First, by identifiying the canonical moments as variables appearing in the CS-decomposition of $X_N$. In the scalar case, this goes back to \cite{sutton2009computing} and \cite{edelman2008beta}. Any effort to generalize this to a block-CS-demposition seems to fail due to non-commutativity of the blocks. The other possible strategy is to follow the path of \cite{Killip1} applying 
the (inverse) Szeg{\H o} mapping. This yields a symmetric measure on the unit circle $\mathbb{T}$. Then apply the Householder algorithm to the corresponding unitary matrix. Unfortunately,  in the matrix case, the Szeg{\H o} mapping does not give a good symmetric measure on $\mathbb{T}$ in the matrix case. We refer to Section 4.2 
 for a discussion of this difficulty. 
In the present paper, we obtain the distribution of canonical moments by directly 
computing the Jacobian of a compound map. The first application maps 
 support points and weights of $\Sigma_N$ to the recursion coefficients.
Then the recursion coefficients are mapped 
 to a suitable Hermitian version of the canonical moments. We give two different ways to compute this distribution. One proof follows by
 direct calculation. The other one is 
 more subtle. 
It uses the relation between the canonical coefficients and the 
matrix Verblunsky coefficients.

This paper is 
organized  as follows.
In Section \ref{sec:matrixmeasures}, we first give notations and explain the different representations for the 
 matrix measures. We also 
 discuss finitely supported matrix measures. 
 In Section \ref{sec:sumrule}, we 
give our new sum rule. 
Section \ref{sec:randommatrices} is devoted to the set up of probability distributions of the matrix models and of the canonical moments. This leads in Section  \ref{sec:largedeviations} to an LDP for the coeffcient side.
 Section \ref{sec:mainproofs} contains the proof of our three main results, subject to technical lemmas, 
whose proofs are postponed to 
 Section \ref{sec:proofs}.

\section{Matrix measures and representation of coefficients}
\label{sec:matrixmeasures}
 
All along this paper, $p$ will be a fixed integer.  A  $p\times p$ matrix measure $\Sigma$ on $\mathbb R$ is a matrix of complex valued Borel measures  on $\mathbb R$ such that for every Borel set $A \subset \mathbb R$ the matrix $\Sigma(A)$ is nonnegative definite, i.e.  $\Sigma (A) \geq 0$.  
When its $k$-th moment is finite, it is denoted by
\[M_k(\Sigma) = \int x^k \dd \Sigma (x), \qquad k \geq 1,\]
writing $M_k$ for $M_k(\Sigma)$ if the measure is clear from context. 
We keep, as much as possible, the notations  close to those of \cite{GaNaRomat}. All matrix measures in this paper will be of size $p\times p$. 
Let ${\bf 1}$  be the $p \times p$ identity  matrix and ${\bf 0}$ be the $p \times p$ null matrix.
For every integer $n$, $I_n$ denotes the $np\times np$ identity matrix. The set of all matrix measures with support in some set $A$ is denoted by $\mathcal{M}_p(A)$, and we let $\mathcal{M}_{p,1}(A) := \{ \Sigma \in \mathcal{M}_p(A):\, \Sigma(A)={\bf 1} \}$ 
denoting the set of normalized measures. 

For the remainder of this section, let $\Sigma\in \mathcal{M}_{p,1}(\mathbb{R})$ 
have compact support. Such a measure $\Sigma$ can be uniquely described by its sequence of moments $(M_1(\Sigma),M_2(\Sigma),\dots)$. Another particular convenient set of parameters characterizing the measure is given by the coefficients in the recursion of orthogonal matrix polynomials, 
introduced in the following subsection. We will follow 
 largely the exposition developped in \cite{damanik2008analytic}. For matrix measures supported by $[0,1]$, there exists, just as in the scalar case, a remarkable decomposition of the recursion coefficients into a set of so-called canonical moments. The parametrization of $\Sigma$ in terms of these canonical moments is one of the main tools 
 for our probabilistic results.

\subsection{Orthogonal matrix polynomials}

The (right) inner product of two matrix polynomials ${\bf f} , {\bf g}$, i.e., polynomials whose coefficients are complex $p\times p$ matrices, is defined by
\[\lal f ,  g \rar = \int f(x)^\dagger d\Sigma(x)\!\ g(x)\,.\]
A matrix measure is called nontrivial, if for any non zero polynomial $P$ we have
\begin{align} \label{nontrivial}
\tr \lal P,P \rar >0 ,
\end{align}
see Lemma 2.1 of \cite{damanik2008analytic} for equivalent characterizations of nontriviality. 
Let us first suppose that $\Sigma$ is nontrivial. Lemma 2.3 of \cite{damanik2008analytic} shows that then $\lal Q,Q\rar $ is positive definite for any monic polynomial $Q$ (with leading coefficient ${\bf 1}$). 
We may then apply the Gram-Schmidt procedure to $\{\idbold, x\idbold, \dots\}$ and obtain a sequence of monic  matrix polynomials $P_n, n \geq 0$, where $P_n$ has degree $n$ and which are orthogonal with respect to $\Sigma$, that is, $ \lal P_n ,  P_m \rar={\bf 0}$ if $m\neq n$.
The polynomials satisfy the recurrence 
\begin{equation}
\label{recursion1}
xP_n = P_{n+1} + P_n u_n + P_{n-1} v_n , \qquad n \geq 0,
\end{equation}
 where, setting 
\begin{align}  
  \gamma_n := \lal P_n, P_n\rar\,, 
\end{align}
 $\gamma_n$ is Hermitian and positive definite, and for $n\geq 1$
\begin{align} \label{hidden1} 
u_n = \gamma_n^{-1} \lal P_n, xP_n\rar ,  \quad
v_n = \gamma_{n-1}^{-1}\gamma_n , 
\end{align}
with 
 $v_0={\bf 0}$. 
This defines a one-to-one correspondence between the sequence $(u_0,v_1,u_2,\dots)$ and the measure $\Sigma$. 

From the matrix coefficients $u_n,v_n$, we can then define a sequence of very useful Hermitian matrices.
We first define matrices related to 
 orthonormal polynomials recursion. Let for $n\geq 0$
\begin{align}
\label{defAAn}
\tilde{\mathcal{A}}_{n+1} &:= \gamma_{n}^{1/2} v_{n+1}\gamma_{n+1}^{-1/2} =\gamma_{n}^{-1/2}\gamma_{n+1}^{1/2} ,  \\
 \ \mathcal{B}_{n+1} &:= \gamma_n^{1/2}u_n \gamma_n^{-1/2} = \gamma_{n}^{-1/2} \lal P_n, x P_n \rar \gamma_n^{-1/2} \label{defBn}.
\end{align}
Obviously, setting
\[{\bf p}_n = P_n \gamma_n^{-1/2}, \qquad n\geq 0,\]
 defines a sequence of matrix orthonormal polynomials. These polynomials 
satisfy
the recursion
\begin{align}
\label{recursion2}
x {\bf p}_n = {\bf p}_{n+1}\tilde{\mathcal A}_{n+1}^\dagger +{\bf p}_n \mathcal B_{n+1} + {\bf p}_{n-1}\tilde{\mathcal A}_{n},  \qquad n\geq 0),
\end{align}
taking ${\bf p}_{-1}={\bf 0}$. The matrices $\tilde{\mathcal{A}}_{n}$ and $\mathcal{B}_n$ play the role of matrix Jacobi coefficients in the following sense. Define the infinite block-tridiagonal matrix 
\begin{align} \label{jacobimatrix}
J = 
 \begin{pmatrix} 
 \mathcal B_1 &  \tilde{\mathcal A}_1    &                  \\
                \tilde{\mathcal A}_1^\dagger & \mathcal B_2    & \ddots           \\
                    & \ddots & \ddots   
\end{pmatrix} .
\end{align}
On the space of matrix polynomials, the map $f \mapsto (x \mapsto xf(x))$  is a right homomorphism, represented in the (right-module) basis ${\bf p}_0, {\bf p}_1, \dots$  by the matrix $J$. 
Moreover, the measure $\Sigma$ is nothing more than the spectral measure of the matrix $J$ defined through its moments by  
\begin{align*}
e_i^* \int x^k\, d\Sigma(x)e_j = e^*_iJe_j, \qquad i,j =1, \dots, p.
\end{align*} 
 (See for example Theorem 2.11 of \cite{damanik2008analytic}). 

The matrix $\mathcal{B}_n$ is Hermitian and we define the Hermitian square of $\tilde{\mathcal{A}}_n$ by
\begin{equation} \label{defAn}
\mathcal A_n = \tilde{\mathcal A_n}\tilde{\mathcal A_n}^\dagger= \gamma_{n-1}^{-1/2} \gamma_n  \gamma_{n-1}^{-1/2}\,.
\end{equation}
Note that $\mathcal A_n$ is Hermitian positive definite.

\subsection{Measures on $[0,1]$}

Now suppose that $\Sigma$ is a nontrivial matrix measure supported by a subset of 
 $[0,1]$. We present two (equivalent) ways to parametrize $\Sigma$, extending 
the corresponding parametrization of the scalar case. The first one uses the canonical moments, the second one uses the Szeg{\H o} mapping and Verblunsky coeffcients.

\subsubsection{Encoding via canonical coefficients}

Dette and Studden \cite{destu02} proved the following matrix version of Favard's Theorem for measures on $[0,1]$: If $\Sigma$ has support in $[0,1]$, there exist matrices $U_n$, $n\geq 1$, such that the recursion coefficients 
defined in \eqref{recursion1} can
may be decomposed as 
\begin{align} \label{decomposition}
u_n =\zeta_{2n+1}+\zeta_{2n}, \qquad v_n =   \zeta_{2n-1}\zeta_{2n} , \qquad n\geq 1,
\end{align}
where $\zeta_0={\bf 0}, \zeta_1=U_1$ and for $n>1$
\begin{align}\label{recursion3}
\zeta_n = ({\bf 1}-U_{n-1})U_n\,. 
\end{align}

Moreover,  $U_n$ has the following geometric interpretation. Suppose $M_1,\dots ,M_{n-1}$ are the first $n-1$ matrix moments of some nontrivial matrix probability measure on $[0,1]$. Then there exist Hermitian matrices $M_n^-$, $M_n^+$, which are upper and lower bounds for the $n$-th matrix moment. 
More precisely, $M_1,\dots ,M_n$ are the first $n$ moments of some nontrivial measure with support in $[0,1]$, if and only if 
\begin{align} \label{momentineq}
M_n^-< M_n < M_n^+ . 
\end{align}
Here we use the partial
Loewner ordering, that is, $A> B$ ($A\geq B$) for Hermitian matrices $A,B$, if and only if $A-B$ is positive (non-negative) definite. 
Then, if $M_n$ are the moments of a nontrivial measure, 
 the following representation holds: 
\begin{align} \label{canonicalmoment}
U_n = (M_n^+ - M_n^-)^{-1} (M_n - M_n^-)\,.
\end{align}
So that, 
 $U_n$ is the relative position of $M_n$ within the set of all possible $n$-th matrix moments, given the matrix moments of lower order. For this reason, $U_n$ is also called \emph{canonical moment}. Let us define
\begin{align} \label{range}
R_n = M_n^+-M_n^-, \qquad H_n = M_n - M_n^-, 
\end{align}
so that $U_n = R_n^{-1} H_n$. 
A Hermitian version of the canonical moments can be defined by
\begin{align}\label{canonicalmoment2}
\mathcal{U}_n = R_n^{1/2} U_n R_n^{-1/2} = R_n^{-1/2} H_n R_n^{-1/2} . 
\end{align}
The matrices $\mathcal{U}_n$ have been considered previously in \cite{dette2012matrix}, to study asymptotics 
in the random matrix moment problem. 
Note that $U_n$ and $\mathcal{U}_n$ are similar and 
\[ {\bf 0} < \mathcal U_n < {\bf 1}\,.\]
Finally, we remark that $M_n^-,M_n^+$ are continuous functions of $M_1,\dots , M_{n-1}$, and that 
\begin{align}
\label{hgamma}
H_{2n} = \gamma_n .
\end{align}

\subsubsection{Encoding via Szeg{\H o} mapping}
\label{susec:Szego}
 
The Szeg{\H o} mapping is two-one from $\mathbb T =  \{z  \in\mathbb C : |z| =1\}$ to $[-2, 2]$ defined by 
\begin{align}
\label{Smap}
z \in  \mathbb T\mapsto z+ z^{-1} \in [-2,2]\,.\end{align} 
 This induces a bijection $\Sigma \mapsto \hbox{Sz}(\Sigma)$ between matrix probability measures on $\mathbb T$ invariant by $z \mapsto \bar z$ and  matrix probability measures  on $[-2, 2]$. 
On $\mathbb T$, a matrix measure is characterized by the system of its matricial Verblunsky coefficents, ruling the recursion of (right) orthogonal polynomials. When the measure is invariant, the
Verblunsky coefficients $(\abold_n)_{n\geq 0}$ are Hermitian   (\cite{damanik2008analytic} Lemma 4.1) and satisfy ${\bf 0} < \abold_n^2 < \idbold$ for every $n$.

The Verblunsky coefficients of such a matrix probability measure on $\mathbb T$  and the Jacobi coefficients  of the corresponding matrix measure on $[-2, 2]$ are connected by the Geronimus relations (\cite{damanik2008analytic} Theorem 4.2). It is more convenient here 
 to consider the matrix measure  on $[0,1]$  denoted by $\widetilde{\operatorname{Sz}}(\Sigma)$, obtained by
pushing forward  Sz$(\Sigma)$  by the affine mapping $x \mapsto (2-x)/4$.

For $n \geq 0$, let 
 $\abold_n$ be the Verblunsky coefficient of $\Sigma$ and $\mathcal{U}_{n+1}$ the Hermitian 
 canonical moment of $\widetilde{\operatorname{Sz}}(\Sigma)$. Then, the following equality 
 holds: 
\begin{align}
\label{DeWag}
\abold_n  = 2 \mathcal U_{n+1} - 1\,.
\end{align}
The correspondance between the two above encodings is proven in \cite{dewag09}, Theorem 4.3, for real-valued matrix measure. The general complex case is considered in \cite{jensdiss}.  

\begin{rem}
In the scalar case, the canonical parameters $U_n$ 
 can be identified in the CS decomposition (see Edelman-Sutton \cite{edelman2008beta}). In the matrix case, this approach does not seem to work, due to the lack of commutativity.
\end{rem}

\subsection{Finitely supported measures}
\label{sec:finitelysupported}

When the support of $\Sigma$ consists of $N= np$ distinct points, then \eqref{nontrivial} cannot be satisfied for all non zero polynomials and $\Sigma$ is not nontrivial. However, if \eqref{nontrivial} is satisfied for all non zero polynomials of degree at most $n-1$, then actually $\lal Q,Q\rar$ is positive definite for all monic polynomials of degree at most $n-1$, see Lemma 2.3 in \cite{damanik2008analytic}.  This implies that we can use 
 the Gram-Schmidt method to define monic orthogonal polynomials up to degree $n$. Further, 
 $\gamma_k = \lal P_k, P_k\rar$ is positive definite for  $k\leq n-1$. Therefore, the orthogonal polynomials 
allow  also to define the recursion coefficients $u_0, \dots, u_{n-1}; v_1, \dots v_{n-1}$. So that, 
 we can construct $\tilde{\mathcal  A}_1, \dots, \tilde{\mathcal  A}_{n-1};  \mathcal B_1, \dots, \mathcal B_n$ as well, with $\tilde{\mathcal{A}}_k$ nonsingular for $k=1, \dots, n-1$.
Let us denote by $J_{n}$ the $np\times np$ Hermitian block matrix of Jacobi coefficients 
\begin{equation}
\label{jacmatrix0}
J_{n}
= \begin{pmatrix} 
 \mathcal B_1 & \tilde A_1    &         &         \\
               \tilde{\mathcal A}_1^\dagger &\mathcal B_2    & \ddots  &         \\
                    & \ddots & \ddots  & \tilde{\mathcal A}_{n-1} \\
                    &        & \tilde{\mathcal A}_{n-1}^\dagger &\mathcal B_n
\end{pmatrix} .
\end{equation} 
Let $\Sigma^{J_n}$ denote the spectral measure of $J_n$, as defined by \eqref{defmatrixspectral}. The same calculation as in the scalar case shows that the first $2n-1$ moments of $\Sigma^{J_n}$ coincide with those of $\Sigma$. Since these matrix moments determine uniquely the recursion coefficients of monic orthogonal polynomials, the entries of the matrix \eqref{jacmatrix0} are then also the recursion coefficients of orthonormal polynomials  
 for $\Sigma^{J_n}$. 
 
Now, suppose that the support points of $\Sigma$ lie in $[0,1]$. The existence of the canonical moments is tackled in the following lemma, proved in 
Section 7.
It requires some additional assumption and is not so obvious. 

\begin{lem} \label{lem:existencecanonical}
Suppose $\Sigma \in \mathcal{M}_{p,1}([0,1])$ is such that $\tr \lal P,P\rar >0$ for all non zero polynomials of degree at most $n-1$. Suppose further 
 $\Sigma(\{0\})=\Sigma(\{1\})={\bf 0}$. Then, the matrices $M_k^-,M_k^+$ for $k\leq 2n-1$ still exist 
 and 
they satisfy $M_k^-<M_k<M_k^+$ for $k\leq 2n-1$. Moreover, the matrices 
\begin{align}
U_k = (M_k^--M_k^+)^{-1}(M_k-M_k^-), \qquad 1\leq k \leq 2n-1,
\end{align}
are related to the recursion coefficients $u_0, \dots, u_{n-1}; v_1, \dots v_{n-1}$ of $\Sigma$ as in  
\eqref{decomposition} and
\eqref{recursion3}.
\end{lem}

Lemma \ref{lem:existencecanonical} implies that we may still define the Hermitian variables $\mathcal U_1, \dots, \mathcal U_{2n-1}$, if the measure $\Sigma$ is sufficiently nontrivial. In conclusion, for any measure satisfying the assumptions of Lemma \ref{lem:nontrivial}, we have a one-to-one correspondence between: 
\begin{itemize}
\item matrix moments $M_1,\dots ,M_{2n-1}$, with $M_k^-<M_k<M_k^+$ for $k=1,\dots ,2n-1$,
\item recursion coefficients $\mathcal{B}_1,\dots ,\mathcal{B}_{n}$ as in \eqref{defBn} and positive definite $\mathcal{A}_1,\dots ,\mathcal{A}_{n-1}$ as in \eqref{defAn},
\item canonical moments $\mathcal{U}_1,\dots ,\mathcal{U}_{2n-1}$ as in \eqref{canonicalmoment2}, with ${\bf 0} < \mathcal{U}_k < {\bf 1}$ for $k\leq 2n-1$. 
\end{itemize}

\section{The Jacobi sum rule}
\label{sec:sumrule}

The reference measure for the sum rule in the Jacobi case is the matricial version of the 
Kesten-McKay law.  
In the scalar case, this measure is defined for parameters $\kappa_1,\kappa_2\geq 0$ by 
\begin{align*}
\hbox{KMK}(\kappa_1,\kappa_2)(dx) = \frac{2+ \kappa_1 + \kappa_2}{2\pi}\frac{\sqrt{(u^+ -x)(x- u^-)}}{x(1-x)} \ \mathbbm{1}_{(u^-, u^+)}(x) dx\, ,
\end{align*}
where 
\begin{align}
\label{upm}
u^\pm := \frac{1}{2} + \frac{\kappa_1^2 - \kappa_2^2 \pm 4 \sqrt{(1+\kappa_1)(1+\kappa_2)(1+\kappa_1+\kappa_2)}}{2(2+\kappa_1+\kappa_2)^2} .
\end{align}
It appears (sometimes in other parametrizations) as a limit law for spectral measures of regular graphs (see \cite{mckay1981expected}), as the asymptotic eigenvalue distribution of the Jacobi ensemble (see \cite{dette2009some}), or in the study of random moment problems (see \cite{dette2018universality}). For $\kappa_1=\kappa_2=0 $, it reduces to the arcsine law. 
The matrix version is then denoted by 
\begin{align}
\Sigma_{\KMK(\kappa_1, \kappa_2)} := \KMK(\kappa_1, \kappa_2)\cdot \idbold .
\end{align}
The canonical moments of $\Sigma_{\KMK(\kappa_1, \kappa_2)}$ of even/odd order are given by 
\begin{align} \label{KMKcanonical}
U_{2k} = U_e := \frac{1}{2+\kappa_1+\kappa_2} \cdot \idbold ,  
\qquad  U_{2k-1} = U_o := \frac{1+\kappa_1}{2+\kappa_1+\kappa_2}\cdot \idbold\,.
\end{align}
(See \cite{gamboa2011large} Sect. 6 
 for the scalar case, which can obviously be extended 
 to the matrix case.) 

Both sides of 
our sum rule (Theorem \ref{LDPSR}) will only be finite for measures satisfying a certain condition on their support, 
related to the Kesten-McKay law. Let $I=[u^-,u^+]$. We define $\mathcal{S}_p = \mathcal{S}_p(u^-,u^+)$
as the set of all bounded nonnegative matrix
measures $\Sigma \in\mathcal{M}_p(\mathbb{R})$ 
that can be written as
\begin{align}\label{muinS0}
\Sigma = \Sigma_{I} +  \sum_{i=1}^{N^+} \Gamma_i^+ \delta_{\lambda_i^+} + \sum_{i=1}^{N^-} \Gamma_i^- \delta_{\lambda_i^-},
\end{align}
where $\operatorname{supp}(\Sigma_I)\subset I$, $N^-,N^+\in\N_0\cup\{\infty\}$, $\Gamma_i^\pm$ are rank $1$ Hermitian matrices 
 and
\begin{align*}
0 \leq\lambda_1^-\leq\lambda_2^-\leq\dots <u^- \quad \text{and} \quad 1 \geq \lambda_1^+\geq\lambda_2^+\geq\dots >u^+\, .
\end{align*}
We assume that 
 $\lambda_j^-$ converges towards $u^-$ (resp. $\lambda_j^+$ converges to $u^+$) whenever $N^-$ (resp. $N^+$) is not finite. 
An atom outside  $[\alpha^-,\alpha^+]$ may appear several times in the decomposition. Its multiplicity is the rank of the total matrix weight that is decomposed in a sum of rank $1$ matrices.
We also define \[\mathcal{S}_{p,1}=\Sr_{p,1}(u^-,u^+):=\{\Sigma \in \mathcal{S}_p(u^-,u^+) |\, \Sigma(\R)={\bf 1}\}\,.\]

Furthermore, the spectral side of the sum rule of Theorem \ref{LDPSR} involves 
the relative entropy with respect to the central measure. If $\Sigma$ has the Lebesgue decomposition
\begin{align}
\label{gamboa}\Sigma(dx) = h(x) \Sigma_{\KMK}(dx) + \Sigma^s (dx) ,
\end{align}
with $h$ positive $p\times p$ Hermitian and $\Sigma^s$ singular with respect to $\Sigma_{\KMK}$, then we 
define the Kullback-Leibler distance of $\Sigma_{\KMK}$ with respect to $\Sigma$ as
\[\mathcal K(\Sigma_{\KMK} \, | \, \Sigma) = - \int \log \det h(x)  \Sigma_{\KMK}(dx)\,.\]
Let us remark that if $K(\Sigma_{\KMK} \, | \, \Sigma)$ is finite, then $h$ is positive definite almost everywhere on $I$, which implies that $\Sigma$ is nontrivial. Conversely, if $\Sigma $ is trivial, then $K(\Sigma_{\KMK} \, | \, \Sigma)$ is infinite.

Finally, for the contribution of the outlying support points, we define two functionals
\begin{align} \label{outlierF+}
{\mathcal F}_J^+(x) = \begin{cases} \ \displaystyle \int_{u^+}^x \frac{\sqrt{(t -  u^+)(t - u^-)}}
{t(1-t)}\!\ dt  & \mbox{ if} \  u^+ \leq  x \leq 1, \\
\ \infty & \mbox{ otherwise.}
\end{cases}
\end{align}
 Similarly, let
\begin{align} \label{outlierF-}
{\mathcal F}_J^-(x) = \begin{cases}\ \displaystyle \int_x^{u^-} \frac{\sqrt{(u^--t)(u^+ -t)}}
{t(1-t)}\!\ dt & \mbox{ if} \  0 \leq x \leq u^-,\\
\ \infty & \mbox{ otherwise.}
\end{cases}
\end{align}
We are now able to formulate our main result consisting in a sum rule for the matrix Jacobi case.

\begin{thm}
\label{LDPSR}
For $\Sigma \in \mathcal S_{p, 1}(u^-,u^+)$ a nontrivial measure with canonical moments $(U_k)_{k\geq 1}$, we have
\begin{align} 
\label{sumrule}
\mathcal K(\Sigma_{\KMK}\, |\, \Sigma) + \sum_{i=1}^{N^+} \mathcal F^+_J ( \lambda_i^+) + \sum_{i=1}^{N^+} \mathcal F_J (\lambda_i^-) & =
\sum_{k=1}^\infty \mathcal H_o ( U_{2k+1}) +  \mathcal H_e (U_{2k})
\end{align} 
where, for  a matrix $U$ satisfying ${\bf 0} \leq  U \leq {\bf 1}$,
\begin{align} \label{Hoddeven}
\begin{split}
\mathcal H_e(U) &:=  - (\log \det U - \log \det U_e) -(1+\kappa_1+\kappa_2) \left(\log \det ({\bf 1}-U) - \log \det ({\bf 1} - U_e)\right) , \\
\mathcal H_o (U) &:= -(1+ \kappa_1) \left(\log \det U - \log \det U_o\right)  -(1+\kappa_2) \left(\log \det({\bf 1}-U)-\log \det ({\bf 1} -U_o)\right) ,
\end{split}
\end{align}
and where both sides may be infinite simultaneously. If $\Sigma \notin \mathcal S_{p, 1}(u^-, u^+)$, the right hand side equals $+\infty$.
\end{thm}

\begin{rem} \label{rem:sumrulearguments}
The arguments on the right hand side of the sume rule are the canonical moments as they appear in the decomposition of recursion coefficients in \eqref{decomposition} and
\eqref{recursion3}. For some applications, it might be more convenient to work with the Hermitian version as defined in \eqref{canonicalmoment2}. Indeed, since $\mathcal{H}_e, \mathcal{H}_o$ are invariant under similarity transforms, the value of the right hand side does not change when the Hermitian canonical moments $\mathcal{U}_k$ are considered. 

We also 
point out that for trivial measures, $U_k$ or $1-U_k$ will be singular for some $k$ and then the right hand side equals $+\infty$ (see also Theorem \ref{thm:LDPcoefficient}). Since in this case the Kullback-Leibler divergence equals $+\infty$ as well, the equality in Theorem \ref{LDPSR} is also true for trivial matrix measures. 
\end{rem}

As in previous papers, an important consequence of this sum rule 
 a system of  equivalent conditions for finiteness of both sides. 
It is  a {\it gem}, as defined by Simon in \cite{simon2} p.19. The following statement is the {\it gem} implied by Theorem \ref{LDPSR}.
We give equivalent conditions on the matrices $\mathcal U_k$ and the spectral measure, which characterize the finiteness of either side in the sum rule identity. The following corollary is the matrix counterpart of Corollary 2.6 in \cite{magicrules}. It follows immediately from Theorem \ref{LDPSR}, since 
\[\mathcal F^\pm_J (u^\pm \pm h) = \frac{2\sqrt{u^+-u^-}}{3u^\pm (1-u^\pm)}h^{3/2} + o(h^{3/2}) \ \ \ (h \rightarrow 0^+)\]
and, for $H$ similar to a Hermitian matrix, 
\begin{align*}
\mathcal H_e(U_e+H) &= \frac{(2+\kappa_1+\kappa_2)^2(\kappa_1+\kappa_2)}{2(1+ \kappa_1+\kappa_2)} \tr H^2 + o(||H||^2), \\
\mathcal H_o (U_o + H) &=  \frac{(2+\kappa_1+\kappa_2)^2(\kappa_2 -\kappa_1)}{2(1+ \kappa_1)(1 +\kappa_2)} \tr H^2  + o(||H||^2),
\end{align*}
as $||H|| \rightarrow 0$, where $||\cdot||$ is any matrix norm.

\begin{cor} \label{semigemL}
Let $\Sigma$ be a nontrivial matrix probability measure on $[0,1]$ with canonical moments $(U_k)_{k\geq 1}$. Then for any $\kappa_1, \kappa_2 \geq 0$,
\begin{align}
\label{zl2}
\sum_{k=1}^\infty\left[\tr(U_{2k-1} - U_o)^2  + \tr(U_{2k} -U_e)^2 \right] < \infty
\end{align}
if and only if the three following conditions hold:
\begin{enumerate}
\item 
  $\Sigma \in \mathcal S_{p,1} (u^-,  u^+)$
\item $\sum_{i=1}^{N^+} (\lambda_i^+ - u^+)^{3/2} + \sum_{i=1}^{N^-} (u^- - \lambda_i^- )^{3/2}  < \infty$ and additionally, if $N^->0$, then $\lambda_1^- > 0$ and if $N^+>0$, then $\lambda_1^+<1$. 
\item Writing the Lebesgue decomposition of $\Sigma$ as in (\ref{gamboa}), then
\begin{align*}
\int_{u^-}^{u^+} \frac{\sqrt{(u^+-x)(x-u^-)}}{x(1-x)} \log \det (h(x)) dx >-\infty .
\end{align*}
\end{enumerate}
\end{cor}

\section{Randomization: Classical random matrix ensembles and their spectral measures}
\label{sec:randommatrices}

To prove the sum rule of Theorem \ref{LDPSR} by our probabilistic method, we start from some 
 random Hermitian  matrix $X_N$ of size $N=np$. The random spectral measure $\Sigma_N$ associated with ($X_N; e_1,\dots ,e_p)$, 
is defined through its matrix moments: 
\begin{align}\label{defspectralmeasure}
 M_k(\Sigma_n)_{i,j} = e_i^\dagger X_N^k e_j ,  \qquad k\geq 0,\ 1 \leq i, j\leq p,
\end{align}
where $e_1,\dots ,e_N$ is the canonical basis of $\mathbb C^N$. From the spectral decomposition of $X_N$, we see that the matrix measure $\Sigma_N$ is 
\begin{align} \label{defspectralmeasure2}
\Sigma_N = \sum_{j=1}^{N} \v_j \v_j^\dagger \delta_{\lambda_j} ,
\end{align}
where the support is given by the eigenvalues of $X_N$ and $\v_j$ is the projection of a unitary eigenvector corresponding to the eigenvalue $\lambda_j$ on the subspace generated by $e_1, \dots, e_p$.
A sum rule is then a consequence of two LDPs for the sequence $(\Sigma_N)_n$, 
the first one when the measure is encoded by its support and the weight, as in \eqref{defspectralmeasure2}, and the second one when the measure is encoded by its recursion coefficients. 
The two following questions are therefore crucial: 
\begin{itemize}
\item What is the joint distribution of $(\lambda_1, \dots, \lambda_N; \v_1, \dots, \v_N)$?
\item What is the distribution of the matricial recursion or canonical coefficients? 
\end{itemize}
The answer to the first question is now classical (see \cite{mehta2004random} or \cite{agz}), when $X_N$ is 
chosen according to a density (the joint density of all real entries, up to symmetry constaint) proportional to 
\begin{align} \label{generalpotential}
\exp \big( -N \tr V(X) ) ,
\end{align}
for some potential $V$. In this case, the eigenvalues follow a log-gas distribution and independently, the eigenvector matrix is Haar distributed on the unitary group. In \cite{GaNaRomat}, the authors considered such general potentials and proved an LDP 
 using the encoding by eigenvalues and weights. 
For $X_N$ distributed according to the Hermite and Laguerre ensemble, it is also possible to answer the second question and derive the  LDPs in both 
encodings. Remarkably, the recursion coefficients in the Hermite case are independent and are $p\times p$ matrices of the Hermite and Laguerre ensemble. In the Laguerre case, Hermitian version of the matrices $\zeta_k$ as in \eqref{decomposition} are Laguerre-distributed. 

In this section, we give the answer to the second question, when $X_N$ is a matrix of the Jacobi ensemble. We first introduce all classical ensembles.

\subsection{The classical ensembles: GUE, LUE, JUE}

We denote by $\mathcal N(0, \sigma^2)$ the centered Gaussian distribution with variance $\sigma^2>0$. 
A random variable $X$ taking values in $\mathcal H_N$, the set of all Hermitian $N \times N$ matrices, is distributed according to the Gaussian unitary ensemble $\GUE_N$, if all real diagonal entries are 
distributed as $\mathcal N(0, 1)$ 
 and the real and imaginary parts of off-diagonal variables are independent and $\mathcal N(0, 1/2)$ distributed (also called complex standard normal distribution). 
All entries are assumed to be independent up to symmetry and conjugation. 
The random matrix $X$ has then a density as in \eqref{generalpotential} with $V(x) = \tfrac{1}{2}x^2$. The joint density of the (real) eigenvalues $\lambda = (\lambda_1,\dots ,\lambda_N)$ of $X$
is 
\begin{align} \label{evg}
g_G(\lambda) = c_r^H \Delta(\lambda)^2  \prod_{i=1}^N e^{- \lambda_i^2/2}.
\end{align}
where 
\[\Delta(\lambda) =  \prod_{1\leq  i < j\leq N} |\lambda_i - \lambda_j|\]
is the Vandermonde determinant.

By analogy with 
 the scalar $\chi^2$ distribution, the Laguerre ensemble is the distribution of the ''square'' of Gaussian matrices. 
More precisely, if $a$ is a nonnegative integer and if  
$G$ denotes a $N \times (N+ a)$ matrix with independent complex standard normal entries, then $X= GG^\dagger$  
is said to be distributed according to the Laguerre ensemble $ \LUE_N(N+ a)$. Its density (on the set $\mathcal H_N^+$ of positive definite Hermitian matrices) is proportional to
\[(\det X)^a \exp \big( -\tfrac{1}{2}\tr\!\ X\big) \,.\]
The eigenvalues density in this case is
\begin{align} \label{evl}
g_L(\lambda) = c_{N,a}^L \Delta(\lambda)^2  \prod_{i=1}^N \lambda_i^{a} e^{-\lambda_i} \mathbbm{1}_{\{ \lambda_i>0\} }.
\end{align}

For $a,b$ nonnegative integers, 
let $L_1$ and $L_2$ be independent matrices distributed according to $\LUE_N (N+a)$ and $\LUE_N (N+ b)$, respectively.
Then the Jacobi ensemble $\JUE_N(a,b)$ is the distribution of 
\begin{align} \label{defjac}
X = ( L_1 + L_2)^{-1/2} L_1 (L_1 + L_2)^{-1/2}  .
\end{align}
Its density on the set  of Hermitian $N\times N$ matrices satisfying $0 < X < I_N$ is proportional to
\begin{align} \label{defjac2}
\det X^a \det (I_N-X)^b . 
\end{align}
The  density of the eigenvalues $(\lambda_1, \dots, \lambda_N)$  is then given by
\begin{align} \label{evj}
g_J(\lambda) = c_{N,a,b}^J | \Delta (\lambda )|^{2} \prod_{i=1}^N \lambda_i^{a} (1-\lambda_i)^{b} \mathbbm{1}_{\{ 0<\lambda_i<1 \} } .
\end{align}
By extension we will say that 
$X$ is distributed according to $\JUE_N(a,b)$, if it has density \eqref{defjac2}, for general real parameters $a,b\geq 0$.

As mentioned above, in all three cases the eigenvector matrix is independent of the eigenvalues and Haar distributed on the group of unitary matrices. As a consequence, the matrix weights in 
 the spectral measure (see \eqref{defspectralmeasure2}) have a distribution which is a matrical generalization of the Dirichlet distribution. Let us denote the distribution of $(\v_1 \v_1^\dagger, \dots ,\v_N \v_N^\dagger )$ by $\mathbb{D}_{N,p}$. It was shown in \cite{FGARop}, that this distribution may be obtained as follows:
Let $z_1,\dots ,z_N$ be random vectors in $\mathbb{C}^p$, with all coordinates independent complex standard normal distributed and set $H =z_1z_1^\dagger + \dots + z_1z_1^\dagger$. Then we have the equality in distribution
\begin{align}
\big( \v_1 \v_1^\dagger, \dots ,\v_N \v_N^\dagger \big) \overset{d}{=} \big( H^{-1/2}z_1 z_1^\dagger H^{-1/2}, \dots ,H^{-1/2}z_N z_N^\dagger H^{-1/2} \big) .
\end{align}   
Using this representation, we can prove the following useful lemma, which shows that although our random spectral measures are finitely supported and thus not nontrivial, it is still possible to define the first recursion coefficients or canonical moments. 

\begin{lem} \label{lem:nontrivial}
Let $N=np$ and $\Sigma_N$ be a random spectral measure as in \eqref{defspectralmeasure2}. We assume that 
 there are almost surely $N$ distinct support points and that the weights are $\mathbb{D}_{N,p}$ distributed and independent of the support points. Then, with probability one, for all nonzero matrix polynomials $P$ of degree at most $n-1$, 
\begin{align*}
\tr \lal P,P\rar >0 .
\end{align*}
\end{lem}

\subsection{Distribution of coefficients}

In the following, let $N=np$. If $\Sigma_N$ is a spectral matrix measure of a matrix $X_N\sim\GUE_N$, then, almost surely, the $N$ support points of $\Sigma_N$ are distinct and none of them equal 0 or 1. By Lemma \ref{lem:nontrivial} and the discussion in Section \ref{sec:finitelysupported}, $\Sigma_N$ may be encoded by its first $2n-1$ coefficients in the polynomial recursion. It is known that then the random matrices $ \mathcal B_1, \dots, \mathcal B_n,\mathcal A_1, \dots, \mathcal A_{n-1}$  are independent and 
\[\mathcal A_k \sim  \LUE_p((N-k)p), \qquad  \mathcal B_k \sim \GUE_p .\]
For the Laguerre ensemble, the spectral measure is supported by $[0, \infty)$ and then a decomposition as in  \eqref{decomposition} still holds, where now Hermitian versions of 
 $\zeta_1,\dots ,\zeta_{2n-1}$ 
 are distributed according to the Laguerre ensemble of dimension $p$ with appropriate parameter. 
These results may be seen in \cite{GaNaRomat}, Lemma 6.1 and 6.2. They are extensions 
of the scalar results of Dumitriu-Edelman \cite{dumede2002} and their proofs are in \cite{GaNaRomas}. Since therein they are formulated in a slightly different way, we clarify the arguments in the Hermite case when we prove Theorem \ref{thm:distributioncanonical} below. It is one of our main results, and shows that in the Jacobi case, the matricial canonical moments are independent and again distributed as matrices of the Jacobi ensemble. 

\begin{thm} \label{thm:distributioncanonical}
Let $\Sigma_N$ be the random  spectral matrix measure associated with the $\JUE_N(a,b)$ distribution. Then,
 the Hermitian canonical moments $\mathcal U_1,\dots ,\mathcal{U}_{2n-1}$ are independent and for $k=1, 2, \dots, n-1$,
\begin{align}
\mathcal U_{2k-1}\sim \JUE_p(p(n-k) + a, p(n-k) + b),\quad  \mathcal U_{2k} \sim  \JUE_p(p(n-k-1), p(n-k) +a+b) 
\end{align}
and $\mathcal U_{2n-1}\sim \JUE_p(a, b)$.
\end{thm}

The Jacobi scalar case was solved by Killip and Nenciu \cite{Killip1}. They used the inverse Szeg{\H o} mapping and actually considered the symmetric random measure on $\mathbb T$  as the spectral measure of $(U; e_1)$ where $U$ is an element of 
 $\mathbb S\mathbb O(2N)$ and $e_1$ is the first vector of the canonical basis. This measure 
may be written as
\[\mu= \sum_{k=1}^N\w_k \left(\delta_{e^{i\theta_k}} +\delta_{e^{-i\theta_k}}\right)\,.\]
Under the Haar measure,
 the support points (or eigenvalues) have 
 the joint density 
 proportional to 
\[\Delta(\cos \theta_1, \dots, \cos \theta_N)^2\]
and the weights are Dirichlet distributed. This induces for the 
pushed forward eigenvalues a density proportional to
\[\Delta(\lambda)^2 \prod_{i=1}^N\lambda_i^{-1/2} (1- \lambda_i)^{-1/2}\]
Then they used a ''{\it magic relation}" to get rid of the factor $\prod \lambda_i^{a-1/2}  (1-\lambda_i)^{b-1/2}$.

If we consider the matricial case, i.e. if we sample $U$ according to the Haar measure on  $\mathbb S\mathbb O(2Np)$ with $(p \geq 2)$, the matrix spectral measure of $(U; e_1, \dots, e_p)$ is now
\[\Sigma =\sum_{k=1}^{N} \left(\w_k  \delta_{e^{i\theta_k}} + \bar\w_k  \delta_{e^{-i\theta_k}}\right) ,\]
the eigenvectors of conjugate eigenvalues being conjugate of each other.
 Unfortunately, this measure is symmetric (i.e. invariant by $z \mapsto \bar z$) only in the scalar case $p=1$, which prohibits the use of the Szeg{\H o} mapping. To find the distribution of the canonical moments, we have to follow another strategy. 
First, we will use
 the explicit relation between the distribution of eigenvalues and weights and the distribution of the recursion coefficients, when sampling the matrix in  the Gaussian ensemble. Then we will 
 compute the Jacobian of the mapping from recursion coefficients to canonical moments using the representation 
in terms of moments as in \eqref{canonicalmoment2}.

\section{Large deviations}
\label{sec:largedeviations}

In order to be self-contained, let us recall the definition of a large deviation principle. For a general reference on large deviation statements we refer to the book 
 \cite{demboz98} or to the Appendix D of \cite{agz}.

Let $E$ be a topological Hausdorff space with Borel $\sigma$-algebra $\mathcal{B}(E)$. We say that a sequence $(P_{n})$ of probability measures on $(E,\mathcal{B}(E))$ satisfies the large deviation principle (LDP) with speed $a_n$ and
rate function $\mathcal{I} : E \rightarrow [0, \infty]$  if:
\begin{itemize}
\item [(i)] $\mathcal I$ is lower semicontinuous.
\item[(ii)] For all closed sets $F \subset E$: $\displaystyle \qquad 
\limsup_{n\rightarrow\infty} \frac{1}{a_n} \log P_{n}(F)\leq -\inf_{x\in F}\mathcal{I}(x) $
\item[(iii)] For all open sets $O \subset E$: $\displaystyle \qquad 
\liminf_{n\rightarrow\infty} \frac{1}{a_n} \log P_{n}(O)\geq -\inf_{x\in O}\mathcal{I}(x) $
\end{itemize}
The rate function $\mathcal{I}$ is good if its level sets
$\{x\in E |\ \mathcal{I}(x)\leq a\}$ are compact for all $a\geq 0$. 
We say that a sequence of $E$-valued random variables satisfies an LDP if their distributions satisfy an LDP. 

It was shown in Theorem 3.2 of \cite{GaNaRomat}, that the sequence of matrix spectral measures $\Sigma_N$ of the Jacobi ensemble $\JUE_N(\kappa_1N,\kappa_2N)$ satisfies an LDP with speed $N$ 
 and good rate function equal to the left hand side of the sum rule in Theorem \ref{LDPSR}. The LDP for the coefficient side is given in the following theorem. Its proof is independent of the one given 
 in \cite{GaNaRomat}.

\begin{thm} \label{thm:LDPcoefficient}
Let $\Sigma_N$ be a random spectral matrix measure of the Jacobi ensemble $\JUE_N(\kappa_1 N, \kappa_2 N)$, with $\kappa_1,\kappa_2\geq 0$ and $N=pn$. Then the sequence $(\Sigma_N)_N$ satisfies the LDP in $\mathcal{M}_{p,1}([0,1])$, with speed $N$ and good rate function 
\begin{align}
\mathcal{I}_J(\Sigma) =\sum_{k=1}^\infty \mathcal H_o (U_{2k-1}) +  \mathcal H_e (U_{2k})
\end{align} 
for nontrivial $\Sigma$, where $\mathcal H_o$ and $\mathcal H_e$ are defined in \eqref{Hoddeven} and $U_k, k\geq 1$ are the canonical moments of $\Sigma$. If $\Sigma$ is trivial, then $\mathcal{I}_J(\Sigma)=+\infty$.
\end{thm}

The following lemma shows an LDP for the Jacobi ensemble of fixed size. It is crucial in proving the LDP for the canonical moments and consequently Theorem \ref{thm:LDPcoefficient}. 

\begin{lem} \label{prop:LDPsingleU}
For $\alpha, \alpha' > 0$ suppose that $X_n \sim \JUE_p (\alpha n + a , \alpha' n + b)$. Then $(X_n)_n$ satisfies the LDP in the set of Hermitian $p\times p$ matrices, with speed $n$ and good rate function $I_{\alpha,\alpha'}$
where
\begin{align}
\label{defIalpha}
I_{\alpha,\alpha'}(X) = - \alpha \log \det X - \alpha' \log \det(\idbold -X) + p\alpha \log \frac{\alpha}{\alpha + \alpha'} + p \alpha' \log \frac{\alpha'}{\alpha + \alpha'}
\end{align}
for $\bold{0}<X<\idbold$ and $I_{\alpha,\alpha'}(X)=\infty$ otherwise.
\end{lem} 

The proof of Lemma \ref{prop:LDPsingleU} makes use of the explicit density and follows as Proposition 6.6 in \cite{GNROPUC}. 

\section{Proof of the main results}
\label{sec:mainproofs}

In this section we prove our three main results in the order of their dependence. First, Theorem \ref{thm:distributioncanonical}
 provides  the distribution of the canonical moments 
 for the Jacobi ensemble, then Theorem \ref{thm:LDPcoefficient} shows the  LDP 
 for the spectral measure of the Jacobi ensemble, and finally  Theorem \ref{LDPSR} 
establishes the sum rule for the Jacobi case. For these three proofs, we 
use the result of all our technical lemmas, whose proofs are postponed 
to Section 7.

\subsection{Proof of Theorem \ref{thm:distributioncanonical}}

The starting point is the spectral measure 
\begin{align} \label{recallspectralmeasure}
\Sigma_N = \sum_{i=1}^N \v_i \v_i^\dagger \delta_{\lambda_i} , 
\end{align}
when the distribution of $(\lambda,\v) = (\lambda_1,\dots ,\lambda_{np},\v_1,\dots \v_{np})$ is the probability measure proportional to 
\begin{align} \label{ev+weights}
\left(\Delta(\lambda)^2 \prod_{i=1}^N \lambda_i^a (1-\lambda_i)^b \mathbbm{1}_{\{0<\lambda_i<1\} } d\lambda_j\right) d\mathbb D_{N,p}(\v) .
\end{align}
We need to calculate the pushforward of this measure under the mapping $(\lambda, \v) \mapsto \mathcal{U}=(\mathcal{U}_1,\dots ,\mathcal{U}_{2n-1})$ to the Hermitian canonical moments. By Lemma \ref{lem:nontrivial} and Lemma \ref{lem:existencecanonical} this is well-defined and the canonical moments satisfy ${\bf 0} <\mathcal{U}_k <{\bf 1}$. The first step will be 
the computation of
 the pushforward under the mapping $(\lambda, \v)\mapsto (\mathcal{A},\mathcal{B})$, when $(\mathcal A, \mathcal B ) := (\mathcal A_1, \dots, \mathcal A_{n-1}, \mathcal B_1, \dots, \mathcal B_n)$ are the Hermitian recursion coefficients as defined in \eqref{defBn} and \eqref{defAn}. This can be done by considering the corresponding change of measure in the Gaussian case, that is, when $\Sigma_N$ is the spectral measure of a $\GUE_N$-distributed matrix with distribution proportional to 
\begin{align} \label{ev+weightsG}
\left(\Delta(\lambda)^2 \prod_{i=1}^N e^{-\tfrac{1}{2}\lambda_i^2} d\lambda_j\right) d\mathbb D_{N,p}(\v) .
\end{align}
As mentioned in Section \ref{sec:randommatrices}, the correspondence in the Gaussian case was investigated in \cite{GaNaRomat}. Lemma 6.1 therein shows that the spectral matrix measure $\Sigma_N$ is also the spectral matrix measure of the block-tridiagonal matrix
\begin{align}
\label{jacmatrix1}
\hat J_{n}
= \begin{pmatrix} 
  D_1 &  C_1    &         &         \\
               C_1 & D_2    & \ddots  &         \\
                    & \ddots & \ddots  & C_{n-1} \\
                    &        & C_{n-1} &D_n
\end{pmatrix} ,
\end{align}
where $C_k,D_k$ are Hermitian and independent, with $D_k\sim \GUE_p$ and $C_k$ is positive definite with $C_k^2\sim \LUE_p(p(n-k))$. This implies that the Hermitian recursion coefficients $\mathcal{B}_k$ and $\mathcal{A}_k$ are given by $\mathcal{B}_k = D_k$ and $\mathcal{A}_k = C_k^2$, respectively. That is, the pushforward of the measure \eqref{ev+weightsG} under the mapping $(\lambda, \v)\mapsto (\mathcal{A},\mathcal{B})$ is the measure proportional to 
\begin{align}
\left(\prod_{k=1}^n \exp \left( - \frac{1}{2} \tr \mathcal B_k^2\right)\  d\mathcal B_k\right) \left(\prod_{k=1}^{n-1} (\det \mathcal A_k)^{p(n-k-1)} \exp \left(- \frac{1}{2} \tr \mathcal A_k\right)\  d\mathcal A_k\right).
\end{align} 
Here and in the following, $dM$ 
denotes the Lebesgue measure in each of the functionally independent real entries of a Hermitian matrix $M$. 
Since 
\[\tr \hat J_n^2 = \sum_{k=1}^n \tr \mathcal B_k^2 + \sum_{k=1}^{n-1} \tr \mathcal A_k = \sum_{j=1}^N\lambda_j^2 , \]
 we conclude
that the pushforward of the measure 
\begin{align}
\label{pushm}
\left(\Delta(\lambda)^2 \prod_{i=1}^N \mathbbm{1}_{\{0<\lambda_i<1\} } d\lambda_i\right) d\mathbb D_{N,p}(\w)
\end{align}
by the mapping $(\lambda, \w) \mapsto (\mathcal A,  \mathcal B)$ is, up to a multiplicative constant, the measure
\begin{align}
\label{lawAB}
\left(\prod_{k=1}^{n-1} (\det \mathcal A_k)^{p(n-k-1)} d\mathcal A_k\right)  \prod_{k=1}^n d\mathcal B_k\,.
\end{align}
Note that an indicator function is omitted in \eqref{lawAB}, (ensuring that
 the spectral measure is supported by $(0,1)$). This indicator function will appear 
 in the condition ${\bf 0} < \mathcal{U}_k < {\bf 1}$, but it does not play a role in the following arguments. 

Now 
 two steps are remaining. First, we need to compute the pushforward of \eqref{lawAB} under the mapping $(\mathcal A , \mathcal B) \mapsto \mathcal U := (\mathcal U_1, \dots, \mathcal U_{2n-1})$. 
Second, to express the prefactor $\prod_{i=1}^{np} \lambda_i^a(1-\lambda_i)^b$ in \eqref{ev+weights} in terms of $\mathcal U$. This is summarized in the two following technical lemmas, whose proofs are in Section \ref{appendix0} and \ref{appendix}, respectively.

\begin{lem}
\label{push2}
The pushforward  of the measure \eqref{lawAB}
by the mapping $(\mathcal{A},\mathcal{B}) \mapsto \mathcal U$ is, up to a multiplicative constant, the measure
\begin{align}
\left(\prod_{k=1}^{n-1} \det((\idbold-\mathcal{U}_{2k-1})\mathcal{U}_{2k-1})^{p(n-k)} d\mathcal U_{2k-1}\right) \left(  
\prod_{k=1}^{n-1} \det(\idbold-\mathcal{U}_{2k})^{p(n-k)} \det(\mathcal{U}_{2k})^{p(n-k)} d \mathcal U_{2k}\right)
\end{align}
\end{lem}

\begin{lem}
\label{crucial}
\begin{align}
\prod_{i=1}^{np} (1 - \lambda_i) = \prod_{k=1}^{2n-1} \det (\idbold - \mathcal U_k) , \qquad  
 \prod_{i=1}^{np} \lambda_i = \prod_{k=1}^n \det \mathcal U_{2k-1} \prod_{k=1}^{n-1} \det (\idbold -\mathcal U_{2k})\,.
\end{align}
\end{lem}

Gathering these results we see that the pushforward of the measure \eqref{ev+weights} by the mapping $(\lambda, \w) \mapsto \mathcal U$ is, again up to a multiplicative constant,
\begin{align}
\prod_{k=1}^{n}  \det (\mathcal{U}_{2k-1})^{p(n-k)+a}\det(\idbold-\mathcal{U}_{2k-1})^{p(n-k)+b} 
\prod_{k=1}^{n-1} \det (\mathcal{U}_{2k})^{p(n-k-1)}\det(\idbold-\mathcal{U}_{2k})^{p(n-k)+a+b} 
\prod_{k=1}^{2n-1} d\mathcal{U}_k .
\end{align}
That is, the canonical moments are independent and 
\[\mathcal{U}_{2k-1}\sim \JUE_p (p(n-k)+a,p(n-k)+b)  , \qquad \mathcal{U}_{2k}\sim \JUE_p (p(n-k-1),p(n-k)+a+b)\,.\]
This ends the proof.
\hfill $ \Box$

\subsection{Proof of Theorem \ref{thm:LDPcoefficient}}

Let $\Sigma_N$ be the spectral measure of a $\JUE_N(\kappa_1N,\kappa_2,N)$ distributed matrix, with $N=np$ and $\kappa_1, \kappa_2 \geq 0$. By Lemma \ref{lem:nontrivial} and Lemma \ref{lem:existencecanonical}, the first $2n-1$ canonical moments $U^{(N)}_k$ $1\leq k\leq 2n-1$ and their Hermitian versions $\mathcal{U}^{(N)}_k$, $1\leq k\leq 2n-1$ are well-defined. They are elements of the space 
\begin{align}
\mathcal{Q}_j= \big\{ (H_1,\dots ,H_{2j-1}) |\, H_j\in\mathcal{H}_p \text{ and } {\bf 0}\leq H_j\leq {\bf 1} \text{ for all } j \big\} .
\end{align}
Let us define the sequence
\begin{align} \label{canonicalsequence}
 \mathcal{U}^{(N)} = \big( \mathcal{U}^{(N)}_{1},\dots ,\mathcal{U}^{(N)}_{2n-1},{\bf 0},\dots \big) ,
\end{align}  
as a random element of 
\begin{align}
\mathcal{Q}_\infty =  \big\{ (H_1,H_2,\dots ) |\, H_j\in\mathcal{H}_p \text{ and } {\bf 0}\leq H_j\leq {\bf 1} \text{ for all } j \big\}, 
\end{align}
which we endow with the product topology. By Theorem \ref{thm:distributioncanonical}, 
\begin{align*}
\mathcal{U}^{(N)}_{2k-1}\sim \JUE_p(p(n-k)+\kappa_1np,p(n-k)+\kappa_2np)
\end{align*}
for $1\leq k\leq n$, and then we apply Lemma \ref{prop:LDPsingleU}, to conclude that the sequence $(\mathcal{U}^{(N)}_{2k-1})_{n}$ satisfies the LDP in $\mathcal{Q}_1$ with speed $n$ and good rate function $I_{p+p\kappa_1,p+p\kappa_2}$. If we instead consider the LDP at speed $N$, the rate function becomes
\begin{align*}
p^{-1}I_{p+p\kappa_1,p+p\kappa_2}(\mathcal{U}) & = - (1+\kappa_1) \log \det (\mathcal{U}) - (1+\kappa_2) \log \det(\idbold -\mathcal{U}) \\
& \quad + p(1+\kappa_1) \log \frac{1+\kappa_1}{2+\kappa_1+\kappa_2} + p (1+\kappa_2) \log \frac{1+\kappa_1}{2+\kappa_1+\kappa_2} ,
\end{align*}
where the right hand side is interpreted as $+\infty$, if we do not have ${\bf 0}< \mathcal{U}< {\bf 1}$. Recalling \eqref{KMKcanonical} and \eqref{Hoddeven}, we see that $p^{-1}I_{p+p\kappa_1,p+p\kappa_2}=\mathcal{H}_o$. 
Turning to the canonical moments of even index, Theorem \ref{thm:distributioncanonical} gives, 
\begin{align*}
\mathcal{U}^{(N)}_{2k}\sim \JUE_p(p(n-k-1),p(n-k)+\kappa_1np+\kappa_2np)
\end{align*}
for $1\leq k\leq n-1$. Then Lemma \ref{prop:LDPsingleU} yields the LDP for $(\mathcal{U}^{(N)}_{2k})_{n}$ in $\mathcal{Q}_1$ with speed $N$ and good rate function $p^{-1}I_{p,p+p\kappa_1+p\kappa_2}$, satisfying
\begin{align*}
p^{-1}I_{p,p+p\kappa_1+p\kappa_2}(\mathcal{U}) & = - \log \det (\mathcal{U}) - (1+\kappa_1+\kappa_2) \log \det(\idbold -\mathcal{U}) \\
& \quad + p \log \frac{1}{2+\kappa_1+\kappa_2} + p (1+\kappa_1+\kappa_2) \log \frac{1+\kappa_1+\kappa_2}{2+\kappa_1+\kappa_2} \\
& = \mathcal{H}_e(\mathcal{U}) .
\end{align*}
Since the canonical moments are independent, we get for any $j\geq 1$, that $( \mathcal{U}^{(N)}_{1},\dots , \mathcal{U}^{(N)}_{2j-1})_{n\geq j}$ satisfies the LDP in $\mathcal{Q}_j$ with speed $N$ and good rate function
\begin{align*}
\mathcal{I}^{(j)}(\mathcal{U}_1,\dots ,\mathcal{U}_{2j-1}) = \mathcal{H}_o(\mathcal{U}_1)+\mathcal{H}_e(\mathcal{U}_2)+\dots +\mathcal{H}_o(\mathcal{U}_{2j-1}) .
\end{align*} 
We can now apply the projective method of the Dawson-G\"artner Theorem (see Theorem 4.6.1 in \cite{demboz98}). It yields the LDP for the full sequence $\mathcal{U}^{(N)}$ in $\mathcal{Q}_\infty$, with speed $N$ and good rate function
\begin{align}
\mathcal{I}_\infty (\mathcal{U}_1,\mathcal{U}_2,\dots ) & = \sup_{j\geq 1} \mathcal{I}^{(j)}(\mathcal{U}_1,\dots ,\mathcal{U}_{2j-1}) 
= \sum_{k=1}^\infty \mathcal{H}_o(\mathcal{U}_{2k-1})+\mathcal{H}_e(\mathcal{U}_{2k}) .
\end{align}
This rate function is finite only if ${\bf 0}<\mathcal{U}_k <{\bf 1}$ for all $k$. In particular, the set where it is finite is a subset of the space
\begin{align}
\widehat{\mathcal{Q}}_\infty = \{ H|\, {\bf 0}< H <{\bf 1} \}^{\mathbb{N}} \cup  \bigcup_{j=1}^\infty \Big(  \{ H|\, {\bf 0}< H <{\bf 1} \}^{2j-1} \times \{ {\bf 0} \}^{\mathbb{N}} \Big).
\end{align}
We also have $\mathcal{U}^{(N)}\in \widehat{\mathcal{Q}}_\infty$ for all $n$, see \eqref{canonicalsequence}. It follows from Lemma 4.1.5 in \cite{demboz98}, that $\mathcal{U}^{(N)}$ also satisfies the LDP in $\widehat{\mathcal{Q}}_\infty$, with speed $N$ and good rate function the restriction of $\mathcal{I}_\infty$ to this space. 

Then, we define the mapping $\psi:\widehat{\mathcal{Q}}_\infty\to \mathcal{M}_{p,1}([0,1])$ as follows. If $\mathcal{U} \in \widehat{\mathcal{Q}}_\infty$ is such that ${\bf 0}<\mathcal{U}_k<{\bf 1}$ for all $k$, there is a unique nontrivial $\Sigma\in \mathcal{M}_{p,1}([0,1])$, such that $\Sigma$ has Hermitian canonical moments $\mathcal{U}$, and we define $\psi(\mathcal{U})=\Sigma$. If $\mathcal{U}$ is such that ${\bf 0}<\mathcal{U}_{2j-1}< {\bf 1}$, but $\mathcal{U}_{k}={\bf 0}$ for $k> 2j-1$, we use the correspondence from Section \ref{sec:finitelysupported}: then there are moments $M_1,\dots ,M_{2j-1}$ with $M_{k}^-<M_k^+$ for $k\leq 2j-1$, and we define $\psi(\mathcal{U})$ as the spectral measure of the block Jacobi matrix $J_j$ as in \eqref{jacmatrix0}, constructed with these moments. That is, $\psi(\mathcal U)$ is the unique spectral measure of such a Jacobi matrix with first canonical moments $\mathcal{U}_1,\dots ,\mathcal{U}_{2j-1}$. Then $\mathcal{U}_n\to \mathcal{U}$ implies that the block-Jacobi matrix of $\psi(\mathcal{U}_n)$ converges entrywise to the block-Jacobi matrix of $\psi(\mathcal{U})$, where the latter one is extended by zeros if $\mathcal{U}$ has less nonzero matricial entries than $\mathcal{U}_n$. This implies that the moments of $\psi(\mathcal{U}_n)$ converge to the moments of $\psi(\mathcal{U})$. Since the convergence of moments of matrix measures on the compact set $[0,1]$ implies weak convergence, the mapping $\psi$ is continuous. 

To 
end the proof, we now apply the contraction principle (Theorem 4.2.1 in \cite{demboz98}). We have $\psi(\mathcal{U}^{(N)}) = \Sigma_N$, and as $\psi$ is continuous, the sequence $(\Sigma_N)_n$ satisfies the LDP in $\mathcal{M}_{p,1}([0,1])$ with speed $N$ and good rate function
\begin{align}
\mathcal{I}_J(\Sigma) = \inf_{\mathcal{U}:\psi(\mathcal{U})=\Sigma} \mathcal{I}_\infty(\mathcal{U}) . 
\end{align}
This infimum is infinite, unless $\Sigma$ is nontrivial, and in this case it is given by $\mathcal{I}_\infty$ evaluated at the unique sequence of canonical moments of $\Sigma$.  
\hfill $\Box$

\subsection{Proof of Theorem \ref{LDPSR}}

Let $\Sigma_N$ be the random spectral matrix measure of a matrix with distribution $\JUE_N(\kappa_1N,\kappa_2N)$, with $\kappa_1,\kappa_2 \geq 0$,
 and suppose $N=np$. This distribution corresponds to a random matrix with potential
\begin{align} \label{Jacobipotential}
V(x) = -\kappa_1\log(x)-\kappa_2 \log(1-x) , 
\end{align}
see \eqref{generalpotential}. In the scalar case $p=1$, the equilibrium measure (the minimizer of the Voiculescu entropy or the limit of $\Sigma_N$) is given by $\KMK (\kappa_1,\kappa_2)$, see \cite{magicrules}, p. 515. For this potential, the assumptions (A1), (A2) and (A3) in \cite{GaNaRomat} are satisfied, with 
matrix equilibrium measure $\Sigma_V=\Sigma_{\KMK (\kappa_1,\kappa_2)}$  and then by Theorem 3.2 of that paper, 
 the sequence $(\Sigma_N)_n$ satisfies the LDP in $\mathcal{M}_{p,1}(\mathbb{R})$ with speed $N$ and good rate function 
\begin{align}
\mathcal{I}_{V}(\Sigma) =  \mathcal K(\Sigma_{\KMK(\kappa_1,\kappa_2)}\!\ |\!\ \Sigma) + \sum_{i=1}^{N^+} \mathcal F^+_V ( \lambda_i^+) + \sum_{i=1}^{N^+} \mathcal F_V^- (\lambda_i^-)
\end{align}
for $\Sigma \in \mathcal{S}_{p,1}(u^-,u^+)$, and $\mathcal{I}_{V}(\Sigma)=+\infty$ otherwise. Here, the functions $\mathcal{F}^\pm_V$ are given by
\begin{align}
\label{rate0}
\mathcal{F}_V^+(x) & = \begin{cases}
\mathcal{J}_V(x) - \inf_{\xi \in \R} \mathcal{J}_V(\xi) & \text{ if } u^+\leq x \leq 1, \\
\infty & \text{ otherwise, } 
\end{cases} \\ \label{rate0b}
\mathcal{F}_V^-(x) & = \begin{cases}
\mathcal{J_V}(x) - \inf_{\xi \in \R} \mathcal{J}_V(\xi) & \text{ if } 1\leq x \leq u^-, \\
\infty & \text{ otherwise, } 
\end{cases}
\end{align}
where $\mathcal{J}_V$ is the effective potential
\begin{align*}
V(x) - 2\int \log |x-\xi|\!\ d \KMK (\kappa_1,\kappa_2)(\xi) . 
\end{align*}
On the one hand, as discussed in Proposition 3.2 of  \cite{magicrules}  (see also the references therein), 
 for $V$ 
 in \eqref{Jacobipotential}, we have $\mathcal{F}_V^\pm=\mathcal{F}_J^\pm$,
(see \eqref{outlierF-} and \eqref{outlierF+}). That is, the rate function $\mathcal{I}_V$ is precisely the left hand side of the sum rule in Theorem \ref{LDPSR}.

On the other hand, as shown in Theorem \ref{thm:LDPcoefficient}, the sequence $(\Sigma_N)_n$ 
 satisfies the LDP with speed $N$ and good rate function $\mathcal{I}_J$. Since a large deviation rate function is unique, we get for any $\Sigma \in \mathcal{M}_{p,1}([0,1])$ the identity
\begin{align*}
\mathcal{I}_V(\Sigma) = \mathcal{I}_J(\Sigma)\,, 
\end{align*}
which is the sum rule of Theorem \eqref{LDPSR}. 
\hfill $\Box$

\section{Proof of the technical lemmas}
\label{sec:proofs}

\subsection{Proof of Lemma \ref{lem:existencecanonical}}

The following statements are true for general nonnegative matrix measures $\Sigma\in \mathcal{M}_p([0,1])$ that are  not necessarily normalized. Let us denote the $n$-th moment space of nonnegative matrix measures on $[0,1]$ 
by 
\begin{align}
\mathfrak{M}_{p,n} = \big\{ (M_0(\Sigma),\dots ,M_n(\Sigma)) |\, \Sigma \in \mathcal{M}_{p}([0,1])\, \big\} \subset \mathcal{H}_p^{n+1}\,.
\end{align}
 A comprehensive study of this 
matrix moment space and the relation between canonical moments and recursion coefficients 
has been addressed in \cite{destu02}. 
Indeed, Theorem 2.7 therein shows that if $(M_0,\dots ,M_{2n-1})$ 
lies in the interior of $\mathfrak{M}_{p,2n-1}$, then the upper and lower bound for $M_k$ satisfies $M_k^-<M_k<M_k^+$ for $1\leq k\leq 2n-1$, and then the canonical moments 
\begin{align}
U_k = (M_k^--M_k^+)^{-1}(M_k-M_k^-), \qquad 1\leq k \leq 2n-1
\end{align} 
are well
 defined. Theorem 4.1 of \cite{destu02} shows that the recursion coefficients  $u_0, \dots, u_{n-1}; v_1, \dots v_{n-1}$ of $\Sigma$ satisfy the decomposition as in  
\eqref{decomposition} and
\eqref{recursion3}. Therefore, the statement of Lemma \ref{lem:existencecanonical} follows once we show that for a measure $\Sigma$ satisfying the assumption of the lemma, $(M_0,\dots ,M_{2n-1})$ is in the interior of the moment space $\mathfrak{M}_{p,2n-1}$. Since this result may be of independent interest, we formulate it as a lemma.

\begin{lem} \label{lem:nontrivialmomentspace}
Let $\Sigma \in \mathcal{M}_{p}([0,1])$ such that 
\begin{align} \label{polynomialnorm}
\tr \lal P,P\rar >0 
\end{align}
for all matrix polynomials $P$ of degree at most $n-1$. Then $(M_0,\dots ,M_{2n-3})$
is in the interior of the moment space $\mathfrak{M}_{p,2n-3}$. If additionally $\Sigma(\{0\})=\Sigma(\{1\})={\bf 0}$, then $(M_0,\dots ,M_{2n-1})$
is in the interior of the moment space $\mathfrak{M}_{p,2n-1}$.
\end{lem}

By the above lemma, there are two sufficient conditions for the existence of the first $2n-1$ canonical moments: either \eqref{polynomialnorm} is satisfied for all polynomials up to degree $n$, or it holds for polynomials up to degree $n-1$ and the additional assumption $\Sigma(\{0,1\})={\bf 0}$ is satisfied. If the condition \eqref{polynomialnorm} fails for some polynomial of degree $n$, then atoms at the boundary can indeed cause the moments to be more ''{extremal}''. This can be made more precise in the scalar case, for which we refer to \cite{DeSt97}, Theorem 1.2.5 and Definition 1.2.10. Suppose $\mu$ is a scalar measure on $[0,1]$ with $n$ support points, then any nonzero polynomial with degree less than $n$ has positive $L^2(\mu)$-norm, but there is a polynomial of degree $n$ with vanishing norm. Then the first $2n-3$ moments will be in the interior of the moment space. On the other hand, the fact that $(M_0,\dots ,M_{2n-1})$ 
lies in the boundary of the moment space is actually equivalent to the fact that $\{0,1\}$ has 
 positive mass. If both 0 and 1 are in the support of $\mu$, then already
 $(M_0,\dots ,M_{2n-2})$ lies at the boundary of the moment space. If exactly one support point is equal to 0 or 1, then the first $2n-2$ moments are interior, but the first $2n-1$ ones are not. If the support contains 0, then $M_{2n-1}=M_{2n-1}^-$, whereas 1 in the support implies $M_{2n-1}=M_{2n-1}^+$. The two versions of $\mu$ are then called the lower and upper principal representation of $(M_0,\dots,M_{2n-2})$, respectively. In the matrix case, the boundary of $\mathfrak{M}_{p,2n-1}$ has a more complicated structure and there is no such equivalence.

\textbf{Proof of Lemma \ref{lem:nontrivialmomentspace}:} 
We again refer to 
 \cite{destu02}. 
 Lemma 2.3 
 says that $(M_0,\dots ,M_{m})$ is an element of $\mathfrak{M}_{p,m}$ if and only if, for all matrices $A_0,\dots, A_m$, such that $Q(x) = A_{m} x^m +\dots + A_0$ is nonnegative definite for all $x\in [0,1]$, we have
\begin{align}\label{inmomentspace}
\tr \sum_{k=0}^{m} A_k M_k  \geq 0 .
\end{align}
Note that the case $A_m={\bf 0}$ is also included. Furthermore, 
$(M_0,\dots ,M_{m})$ is an interior point of $\mathfrak{M}_{p,m}$ if and only if,  for all $A_0,\dots ,A_m$ for which such $Q$ is nonnegative definite on $[0,1]$ and nonzero, we have
\begin{align}\label{inmomentspace2}
\tr \sum_{k=0}^{m} A_k M_k  > 0 .
\end{align}
Theorem 2.5 of \cite{destu02} shows that if the degree of $Q$ is even, say $2\ell$, then such a polynomial can be written as
\begin{align} \label{polynomialdecomp}
Q(x) =  B_1(x)B_1(x)^\dagger + x(1-x)B_2(x)B_2(x)^\dagger , 
\end{align}
where $B_1$ and $B_2$ are matrix polynomials of degree $\ell$ and $\ell-1$, respectively. If the degree of $Q$ is equal to $2\ell-1$, then
\begin{align} \label{polynomialdecomp2}
Q(x) =  xB_1(x)B_1(x)^\dagger + (1-x)B_2(x)B_2(x)^\dagger , 
\end{align}
with $B_1,B_2$ of degree $\ell-1$. 
Let $\Sigma\in \mathcal{M}_p([0,1])$ with $M_k$ the $k$-th moment of $\Sigma$. If $m=2\ell$ and $A_m\neq {\bf 0}$, then, using the decomposition \eqref{polynomialdecomp},
\begin{align} \label{polynomialpositive}
\tr \sum_{k=0}^m A_k M_k & = \tr \int Q(x) d\Sigma(x) \notag \\ 
& = \tr \int B_1(x)B_1(x)^\dagger d\Sigma(x) + \tr \int x(1-x)B_2(x)B_2(x)^\dagger d\Sigma(x)  \notag \\
& = \tr \int B(x)_1^\dagger d\Sigma(x)B_1(x) + \tr \int x(1-x)B_2(x)^\dagger d\Sigma(x)B_2(x)  \notag \\ 
& = \tr\, \lal B_1, B_1\rar + \tr\, \lal p B_2 , B_2\rar ,
\end{align}
where $p(x) = x(1-x)$.

A similar calculation can be made if $m=2\ell-1$ and $A_m\neq {\bf 0}$. Together with the characterizations of the moment space by \eqref{inmomentspace} and \eqref{inmomentspace2}, this implies that for any $\Sigma \in \mathcal{M}_p([0,1])$ and matrix polynomial $B$ 
\begin{align} \label{polcriterion}
\tr\, \lal qB, B\rar \geq 0, 
\end{align}
when $q(x)$ is the scalar polynomial 
 $1,x,1-x$ or $x(1-x)$. 
Furthermore, the first $2m-1$ moments of $\Sigma$ are in the interior of the moment space $\mathfrak{M}_{p,2m-1}$, if 
\begin{align} \label{polcriterion2}
\tr\, \lal qB, B\rar > 0, 
\end{align}
whenever $B$ is nonzero and such that the degree of $q(x)B(x)B(x)^\dagger$ is at most $2m-1$. We remark that this is actually equivalent to the criterion given in \cite{destu02} and stated in terms of Hankel matrices. 

Now suppose that $\Sigma $ is such that $\tr \lal P,P\rar >0$ for all nonzero polynomials $P$ of degree at most $n-1$. We show that then \eqref{polcriterion2} is satisfied whenever the degree of $qBB^\dagger$ 
 is at most $2n-3$. For $q(x)=1$ this is trivially true. In the other cases, 
\begin{align} \label{polynomialpositive2}
\tr\, \lal qB , B\rar = \tr\, \lal qB ,qB\rar + \tr\, \lal qB,(1-q)B\rar . 
\end{align}
Since $qB$ has degree at most $n-1$, the first inner product on the right hand side of \eqref{polynomialpositive2} is 
 positive by assumption. The second one is nonnegative by \eqref{inmomentspace}, since  $q(1-q)BB^\dagger$
 is nonnegative definite on $[0,1]$. This proves that $(M_0,\dots ,M_{2n-3})$ is in the interior of $\mathfrak{M}_{p,2m-1}$. 

Now assume 
 $\Sigma(\{0,1\})={\bf 0}$, we show that then \eqref{polcriterion2} is satisfied whenever $qBB^\dagger$ 
 has degree at most $2n-1$ and $B$ is nonzero. In this case, $B$ is of degree at most $n-1$, and $\tr \lal B, B\rar$ 
 is 
 positive. Using that $\Sigma$ has no mass at $0,1$, 
\begin{align}
\tr\, \lal B ,B\rar = \lim_{\varepsilon \to 0}\, \tr \int_{\varepsilon}^{1-\varepsilon} B(x)^\dagger d\Sigma(x) B(x) ,
\end{align}
and then there exists a $\varepsilon>0$, such that the integral on the right hand side is positive. Since $q(x)\geq \varepsilon(1-\varepsilon)$ on $[\varepsilon,1-\varepsilon]$, and $\int_A B(x)^\dagger d\Sigma B(x)$ is always nonnegative definite, 
\begin{align}
\tr\, \lal qB, B\rar \geq \tr \int_{\varepsilon}^{1-\varepsilon} q(x)B(x)^\dagger d\Sigma(x) B(x)\geq 
\varepsilon(1-\varepsilon)\, \tr \int_{\varepsilon}^{1-\varepsilon} B(x)^\dagger d\Sigma(x) B(x) , 
\end{align} 
which gives a 
 positive lower bound. 
\hfill $ \Box $

\subsection{Proof of Lemma \ref{lem:nontrivial}}

Let us begin by noting that if $z_1,\dots ,z_N$ are random vectors in $\mathbb{C}^p$, independent and complex standard normal distributed, then almost surely, any $p$ of these vectors span $\mathbb{C}^p$. This implies that almost surely, $H=z_1 z_1^\dagger+ \dots + z_N z_N^\dagger$ has full rank. Consider such a realization and let 
\begin{align*}
P(x)=C_{n-1}x^{n-1} + \dots + C_1 x + C_0  
\end{align*}
be a matrix polynomial of degree at most $n-1$. We have 
\begin{align*}
\tr \lal P,P\rar = \tr \sum_{i=1}^N P(\lambda_i)^\dagger \v_i\v_i^\dagger P(\lambda_i) = \sum_{i=1}^N \v_i^\dagger P(\lambda_i)^\dagger P(\lambda_i) \v_i = \sum_{i=1}^N ||P(\lambda_i)\v_i||^2. 
\end{align*}
Suppose that $\tr \lal P,P\rar=0$, then the above calculation shows that for all $i$, $\v_i$ is in the kernel of $P(\lambda_i)$. 
We may rewrite this in matrix form by saying that
\begin{align} \label{nontrivial}
\mathbf{W} \mathbf{P} = \mathbf{0} , 
\end{align}
where $\mathbf{P}$ is $np\times p$ with $\mathbf{P}^\dagger = (C_0,\dots ,C_{n-1})$, and $\mathbf{W}$ is $np\times np$ with
\begin{align*}
\mathbf{W} = 
\begin{pmatrix} 
 \v_1^\dagger & \v_1^\dagger \lambda_1    &  \cdots       & \v_1^\dagger\lambda_1^{n-1}        \\
         \v_2^\dagger & \v_2^\dagger \lambda_2   & \cdots  &     \v_2^\dagger \lambda_2^{n-1}    \\
                \vdots    & \vdots & \ddots  & \vdots \\
         \v_{np}^\dagger & \v_{np}^\dagger \lambda_{np} & \cdots & \v_{np}^\dagger \lambda_{np}^{n-1} 
\end{pmatrix} .
\end{align*}
Now, we show that $\mathbf{W}$ is nonsingular, so that the only solution to \eqref{nontrivial} is $\mathbf{P}=0$, that is, $P$ is the zero polynomial. 
Let $\mathbf{H}$ be the $np\times np$ block-diagonal matrix with blocks $H^{1/2}$ on the diagonal, then $\mathbf{H}$ is nonsingular. The matrix  
$\mathbf{Z}=\mathbf{W}\mathbf{H}$ 
has the same structure as $\mathbf{W}$, except that $\v_i$ is replaced by $z_i$.  
We use an argument similar to what has been done  
 in the proof of Lemma 2.2 in \cite{FGARop}. Conditionally on the eigenvalues, the determinant of $\mathbf{Z}$ is a polynomial in the $np^2$ entries of $z_1,\dots ,z_{np}$. Since they are all independent standard Gaussians, they have a joint density and then either $\det(\mathbf{Z})$ is 0 with probability 0 or it is the zero polynomial. Let us fix $z_{kp+i}=e_i$ for $k=0,\dots n-1$, $i=1,\dots,p$, where $e_1,\dots ,e_p$ is the 
canonical basis of $\mathbb{C}^p$. In this case,
\begin{align}
\mathbf{Z} = 
\begin{pmatrix} 
 e_1^\dagger & e_1^\dagger \lambda_1    &  \cdots       & e_1^\dagger\lambda_1^{n-1}        \\
         e_2^\dagger & e_2^\dagger \lambda_2   & \cdots  &     e_2^\dagger \lambda_2^{n-1}    \\
                \vdots    & \vdots & \ddots  & \vdots \\
                e_p^\dagger & e_p^\dagger \lambda_p & \cdots & e_p^\dagger \lambda_p^{n-1} \\
                e_1^\dagger & e_1^\dagger \lambda_{p+1} & \cdots & e_1^\dagger \lambda_{p+1}^{n-1} \\
                \vdots    & \vdots & \ddots  & \vdots \\
         e_{p}^\dagger & e_{p}^\dagger \lambda_{np} & \cdots & e_{p}^\dagger \lambda_{np}^{n-1} 
\end{pmatrix} .
\end{align}
By reordering rows and columns, this matrix may be transformed into the block diagonal matrix $\widetilde{\mathbf{Z}}$ with $n\times n$  Vandermonde-blocks,
\begin{align}
\widetilde{\mathbf{Z}} = 
\begin{pmatrix}
 1  & \lambda_1 & \cdots & \lambda_1^{n-1}    &  & & & & \\
 1 & \lambda_{p+1} & \cdots & \lambda_{p+1}^{n-1}   &   & & & &  \\
 \vdots &  \vdots &  & \vdots &   & & & &  \\
 1 & \lambda_{(n-1)p+1} & \cdots & \lambda_{(n-1)p+1}^{n-1}    &   & & & &  \\
 &  & &  &  \ddots & & & &  \\
& & & & & 1&  \lambda_p & \cdots & \lambda_p^{n-1} \\
& & & & & 1&  \lambda_{2p} & \cdots & \lambda_{2p}^{n-1} \\
& & & & & \vdots  & \vdots  &  & \vdots \\
& & & & & 1 & \lambda_{np} & \cdots & \lambda_{np}^{n-1} \\
 \end{pmatrix} ,
\end{align} 
which has determinant
\begin{align}
\det(\widetilde{\mathbf{Z}}) = \prod_{k=1}^p \prod_{1<j<n} (\lambda_{jp+k}- \lambda_{ip+k} ) . 
\end{align}
Since the $\lambda_i$ are almost surely disjoint, the matrix $\widetilde{\mathbf{Z}}$ is almost surely non-singular, which implies that $\mathbf{W}$ is almost surely nonsingular. 
\hfill $ \Box$

\subsection{Proof of Lemma \ref{push2}}
\label{appendix0}

We have to compute the Jacobian determinant of the mapping $(\mathcal A, \mathcal B) \mapsto \mathcal U$.
We will do this by using the moments as intermediate variables.
Let us begin by noting that $u_{n-1},U_{2n-1}$ depend on $M_1,\dots ,M_{2n-1}$, but not on any higher moments and $v_n,U_{2n}$ depend only on $M_1,\dots ,M_{2n}$. Since for the similarity transforms in Section \ref{sec:matrixmeasures} we used only matrices depending on moments of strictly lower order, the same statements can by made for the Hermitian versions, where $\mathcal{B}_n,\mathcal{U}_{2n-1}$ depend on $M_1,\dots,M_{2n-1}$ and $\mathcal{A}_n,\mathcal{U}_{2n}$ depend on $M_1,\dots ,M_{2n}$.
We have in particular
\begin{align}\label{diffmoments1}
\frac{\partial(\mathcal{B},\mathcal{A})}{\partial M}:=\frac{\partial (\mathcal{B}_1,\mathcal{A}_1,\dots ,\mathcal{B}_{n})}{\partial(M_1,\dots ,M_{2n-1})}
= \frac{\partial \mathcal{B}_1}{\partial M_1}\times \frac{\partial \mathcal{A}_1}{\partial M_2}\times \cdots \times \frac{\partial \mathcal{B}_{n}}{\partial M_{2n-1}} .
\end{align}
Here, we denote by $\frac{\partial F(M)}{\partial M}$ the Jacobian determinant of the mapping $F:\mathcal{H}_p\to \mathcal{H}_p$, seen as a mapping of all the $p^2$ functionally independent real entries of a matrix in $\mathcal{H}_p$, and with the straightforward generalization to mappings with several such matricial coordinates, see \cite{jacobian}. In particular, Theorem 3.5 in \cite{jacobian} shows that for nonsingular $A$,
\begin{align} \label{cnetraljacobian}
\frac{\partial (AMA^\dagger)}{\partial M} = \det(A)^{2p} . 
\end{align}
Recall that
\begin{align}
\label{recall}
\mathcal A_k = H_{2k-2}^{-1/2}H_{2k}H_{2k-2}^{-1/2}  , \qquad  \mathcal{B}_{k} = H_{2k-2}^{-1/2}\lal x P_{k-1}, P_{k-1}\rar H_{2k-2}^{-1/2}
\end{align}
and  $H_{2k} = M_{2k}-M_{2k}^-$ depends only on $M_1, \dots, M_{2k}$ (see \eqref{hgamma}). Then by \eqref{cnetraljacobian}, 
\begin{align} \label{diffA}
\frac{\partial \mathcal{A}_k}{\partial M_{2k}} 
 = \det( H_{2k-2})^{-p} \ , \ \frac{\partial \mathcal{B}_{k}}{\partial M_{2k-1}} 
 = \det (H_{2k-2})^{-p}\,. 
\end{align}
Putting these together, we get that \eqref{diffmoments1} is given by  
\begin{align}\label{diffmoments1result}
\frac{\partial(\mathcal{B},\mathcal{A})}{\partial M} = \det(H_{2n-2})^{-p}\prod_{k=1}^{n-1} \det(H_{2k-2})^{-2p} . 
\end{align}
To end this first step,
we need to evaluate 
\begin{align}\label{diffmoments3}
\frac{\partial \mathcal{U}}{\partial M} := \frac{\partial (\mathcal{U}_1,\dots ,\mathcal{U}_{2n-1})}{\partial(M_1,\dots ,M_{2n-1})}
= \frac{\partial \mathcal{U}_1}{\partial M_1}\times \cdots \times \frac{\partial \mathcal{U}_{2n-1}}{\partial M_{2n-1}} .
\end{align}
where we have by \eqref{canonicalmoment2}
\begin{align}\label{diffU}
\frac{\partial \mathcal{U}_k}{\partial M_{k}} = \frac{\partial \left(R_{k}^{-1/2} H_{k}R_{k}^{-1/2}\right) }{\partial M_{k}} = \det( R_{k})^{-p} 
\end{align}
and then
\begin{align}\label{diffmoments3result}
\frac{\partial \mathcal{U}}{\partial M} = \prod_{k=1}^{2n-1} \det (R_k)^{-p} \,.
\end{align}
Putting together 
\eqref{diffmoments1result} and \eqref{diffmoments3result}, 
 we have shown that
\begin{align}
\label{interm}
\frac{\partial (\mathcal{B},\mathcal{A})}{\partial \mathcal{U}} =  \det(H_{2n-2})^{-p} \prod_{k=1}^{n-1} \det(H_{2k-2})^{-2p} \prod_{k=1}^{2n-1} \det (R_k)^{p} . 
\end{align}
To express this in terms of the canonical moments, we use
\[R_k=R_{k-1}(\idbold -U_{k-1})U_{k-1} , \qquad  H_{k}=R_{k}U_{k}\,,\]
(see \cite{destu02} formulas (2.19) and (2.16)). 
 Taking determinants, we obtain 
\begin{align}
\label{taking}
\det R_k = \prod_{j=1}^{k-1} \det (\idbold -\mathcal U_{j}) \det \mathcal U_{j} , \qquad \det H_{2k-2}= \det R_{2k-2} \det \mathcal U_{2k-2}\,.
\end{align}
We gather \eqref{interm} and \eqref{taking}, to obtain that the pushforward of the measure \eqref{lawAB} by the mapping $(\mathcal{A}, \mathcal{B}) \mapsto \mathcal U$ has, up to a multiplicative constant, 
the density
\begin{align} \label{jacobianABtoU}
 \prod_{k=1}^{n-1} &\det(\mathcal{A}_k)^{p(n-k-1)} 
 \det(H_{2n-2})^{-p} \prod_{k=1}^{n-1} \det(H_{2k-2})^{-2p} \prod_{k=1}^{2n-1} \det (R_k)^{p} \notag \\
&  = \prod_{k=1}^{n-1} \det(H_{2k-2})^{-p(n-k-1)} \det(H_{2k})^{p(n-k-1)}
 \det(H_{2n-2})^{-p} \prod_{k=1}^{n-1} \det(H_{2k-2})^{-2p} \prod_{k=1}^{2n-1} \det (R_k)^{p} \notag \\
  &=  \prod_{k=1}^{n-1} \det(H_{2k-2})^{-p(n-k)} \det(H_{2k})^{p(n-k-1)}
  \prod_{k=1}^{n} \det(H_{2k-2})^{-p} \prod_{k=1}^{2n-1} \det (R_k)^{p}\notag \\
&= \prod_{k=1}^{n} \det(H_{2k-2})^{-p} \prod_{k=1}^{2n-1} \det (R_k)^{p} , 
\end{align}
where for the second line we used \eqref{recall}, 
and then observe the telescopic product of the determinants of $H_k$.

It remains to express \eqref{jacobianABtoU} in terms of the canonical moments. It's time to use \eqref{taking} to get
\begin{align} \label{jacobianABtoU2}
\prod_{k=1}^{n} \det(H_{2k-2})^{-p} \prod_{k=1}^{2n-1} \det (R_k)^{p} & = \prod_{k=1}^{n-1} \det(R_{2k})^{-p}\det(\mathcal{U}_{2k})^{-p} \prod_{k=1}^{2n-1} \det (R_k)^{p} \notag \\ 
& = \prod_{k=1}^{n-1} \det(\mathcal{U}_{2k})^{-p} \prod_{k=1}^{n} \det(R_{2k-1})^{p} \notag \\ 
& = \prod_{k=1}^{n-1} \det(\mathcal{U}_{2k})^{-p} \prod_{k=1}^{n} \prod_{i=1}^{2k-2} \det(\idbold-\mathcal{U}_i)^p \det (\mathcal{U}_i)^p \notag \\
& = \prod_{k=1}^{n-1} \det(\mathcal{U}_{2k})^{-p} \prod_{k=1}^{n-1} \det((\idbold -\mathcal{U}_{2k-1})\mathcal{U}_{2k-1}(\idbold -\mathcal{U}_{2k})\mathcal{U}_{2k})^{p(n-k)}  \notag \\
& = \prod_{k=1}^{n-1} \det((\idbold -\mathcal{U}_{2k-1})\mathcal{U}_{2k-1})^{p(n-k)}  
\prod_{k=1}^{n-1} \det(-\mathcal{U}_{2k})^{p(n-k)} \det(\mathcal{U}_{2k})^{p(n-k)} .
\end{align}
This ends the proof of Lemma \ref{push2}.
\hfill $ \Box$

\subsection{Two proofs of Lemma \ref{crucial}}
\label{appendix}

It follows from Lemma 2.1 of Duran, Lopez-Rodriguez \cite{duran1996orthogonal}, that the eigenvalues of $J_n$ are precisely the zeros of the $n$-th polynomial orthogonal with respect to $\Sigma$. The quadrature formula of Sinap, van Assche \cite{sinap1996orthogonal} implies that the zeros of this polynomial are equal to the support of the spectral measure. As a consequence,  
\begin{align}
\label{detprod}
\det (J_n) = \prod_{i=1}^{np} \lambda_i , \qquad  \det (I-J_n)=\prod_{i=1}^{np} (1-\lambda_i).
\end{align}
In view of (\ref{detprod}) we have to prove that
\begin{align}  \label{finalchange}
\det (I_n - J_n)  = \prod_{k=1}^{2n-1} \det (\idbold - \mathcal U_k) , \qquad \det J_n =  \left(\prod_{k=1}^n \det \mathcal U_{2k-1}\right) \left(\prod_{k=1}^{n-1} \det (\idbold -\mathcal U_{2k})\right).
\end{align}
We give two proofs. The first one is matricial, using a recursion of Schur complements and the second one is based on the Szeg{\H o} mapping and matrix polynomials on the unit circle.

\subsubsection{First proof}

Using the Schur complement formula (see Theorem 1.1 in \cite{hornbasic}),
\begin{align} \label{expansion}
\det(I_n-J_n) & = \det(I_{n-1}-J_{n-1})\det\big(\idbold-\mathcal{B}_n-(0,\dots ,0,\tilde{\mathcal{A}}_{n-1}^\dagger)(I_{n-1}-J_{n-1})^{-1} (0,\dots ,\tilde{\mathcal{A}}_{n-1}^\dagger)^\dagger \big) \notag \\
& = \det(I_{n-1}-J_{n-1}) \det\left( \idbold -\mathcal{B}_n-\tilde{\mathcal{A}}_{n-1}^\dagger [(I_{n-1}-J_{n-1})^{-1}]_{n-1,n-1}\tilde{\mathcal{A}}_{n-1}\right) \notag \\
& = \det(I_{n-1}-J_{n-1}) \det(\varphi_n) ,
\end{align}
where we wrote $[A]_{i,j}$ for the $p\times p$ sub-block in position $i,j$ and we define
\begin{align}
\varphi_n & = \gamma_n^{-1/2}\left(\idbold-\mathcal{B}_n-\tilde{\mathcal{A}}_{n-1}^\dagger [(I_{n-1}-J_{n-1})^{-1}]_{n-1,n-1}\tilde{\mathcal{A}}_{n-1}\right)\gamma_n^{1/2} \notag \\
& = \gamma_n^{-1/2}\left(\idbold -\gamma_n^{1/2}u_{n-1}\gamma_n^{-1/2}-\gamma_{n}^{1/2}\gamma_{n-1}^{-1/2} [(I_{n-1}-J_{n-1})^{-1}]_{n-1,n-1}\gamma_{n-1}^{-1/2}\gamma_n^{1/2}\right)\gamma_n^{1/2} \notag \\
& = \left( \idbold -u_{n-1}-\gamma_{n-1}^{-1/2} [(I_{n-1}-J_{n-1})^{-1}]_{n-1,n-1}\gamma_{n-1}^{-1/2}\gamma_n\right) \notag \\
& = \left(\idbold -u_{n-1}-\gamma_{n-1}^{-1/2} [(I_{n-1}-J_{n-1})^{-1}]_{n-1,n-1}\gamma_{n-1}^{1/2}v_{n-1}\right) .
\end{align}
Recall the non-Hermitian recursion coefficients $u_n, v_n$ have been defined in \eqref{recursion1} and \eqref{hidden1}. Using again the formula of Schur complements (see Theorem 1.2 in \cite{hornbasic}),
\begin{align}
[(I_{n-1}-J_{n-1})^{-1}]_{n-1,n-1} & = \left(\idbold -\mathcal{B}_{n-1}-(0,\dots ,0,\tilde{\mathcal{A}}_{n-2}^\dagger)(I_{n-2}-J_{n-2})^{-1} (0,\dots ,\tilde{\mathcal{A}}_{n-2}^\dagger)^\dagger\right)^{-1} \notag \\
& =  \left(\idbold -\mathcal{B}_{n-1}-\tilde{\mathcal{A}}_{n-2}^\dagger [(I_{n-2}-J_{n-2})^{-1}]_{n-2,n-2}\tilde{\mathcal{A}}_{n-2}\right)^{-1}\notag\\
&= \gamma_{n-1}^{1/2} \varphi_{n-1}^{-1} \gamma_{n-1}^{-1/2} .
\end{align}
We see that $\varphi_n$ satisfies 
 the recursion
\begin{align}\label{recursion}
\varphi_1= \idbold -u_0 , \qquad  \varphi_n  = \idbold -u_{n-1}-\varphi_{n-1}^{-1} v_{n-1}, \qquad  n \geq 2.
\end{align}
 Let us write $V_k=\idbold -U_k$. Then we claim that the solution to this recursion is given by 
\begin{align} \label{recursionsolution}
\varphi_n 
 = V_{2n-2}V_{2n-1} .
 \end{align}
We prove \eqref{recursionsolution} by induction. For $n=1$, we have by \eqref{decomposition} 
\begin{align}
\varphi_1 = \idbold -  u_0 = \idbold -\zeta_1 = \idbold -U_1=V_1 ,
\end{align}
which agrees with \eqref{recursionsolution} since $V_0=\idbold$. Then,
\begin{align}
\varphi_{n+1} & = \idbold -u_{n}-\varphi_n^{-1} v_n \notag \\
& = \varphi_n^{-1}\left[ \varphi_n -\varphi_n (\zeta_{2n}+\zeta_{2n+1}) - \zeta_{2n-1}\zeta_{2n}\right] \notag \\
& = \varphi_n^{-1}\left[  V_{2n-2}V_{2n-1} - V_{2n-2}V_{2n-1} (V_{2n-1}U_{2n}+V_{2n}U_{2n+1}) - V_{2n-2}U_{2n-1}V_{2n-1}U_{2n}\right] \notag \\
& =  \varphi_n^{-1}V_{2n-2} \left[  V_{2n-1} - V_{2n-1}^2U_{2n}-V_{2n-1}V_{2n}U_{2n+1} - U_{2n-1}V_{2n-1}U_{2n}\right] .
\end{align}
In the last line, we write $U_{2n-1}V_{2n-1}U_{2n}=V_{2n-1}U_{2n}-V_{2n-1}^2U_{2n}$ for the last term,  which then cancels the second term in the brackets and leads to
\begin{align}
\varphi_{n+1} 
& =  \varphi_n^{-1}V_{2n-2} \left[  V_{2n-1} -V_{2n-1}V_{2n}U_{2n+1} - V_{2n-1}U_{2n}\right] \notag \\
& = \varphi_n^{-1}V_{2n-2} V_{2n-1}\left[ \idbold  -V_{2n}U_{2n+1} - U_{2n}\right] \notag \\
& = \varphi_n^{-1}\varphi_n \left[ V_{2n} - V_{2n}U_{2n+1}\right] \notag \\
& = V_{2n}V_{2n+1}.
\end{align}
This proves \eqref{recursionsolution}. We may then calculate recursively for  \eqref{expansion}
\[\det(I_n-J_n)  = \det(I_{n-1}-J_{n-1})\det \varphi_n = \det(\varphi_1\dots \varphi_n) , \]
so that
\begin{align} \label{expansionsolution}
\det(I_n-J_n) 
= \prod_{k=1}^{2n-1}\det V_k = \prod_{k=1}^{2n-1}\det ({\bf 1}-U_k) = \prod_{k=1}^{2n-1}\det ({\bf 1}-\mathcal U_k)  \,.
\end{align} 
For the computation of $
\det J_n$, we make use of a decomposition proven in Lemma 2.1 of \cite{GaNaRomat}. There exists a block bi-diagonal matrix $Z_n$, such that $J_n = Z_nZ_n^\dagger$ and (see the proof in \cite{GaNaRomas}), the block $D_k$ in position $k,k$ of $Z_n$ satisfies 
\begin{align}
\det(D_k) = \det(\zeta_{2n-1})^{1/2} .
\end{align}
Then, this implies 
 \begin{align}
 \det J_n &= (\det Z_n)^2 = \prod_{k=1}^n (\det D_k)^2 
= \prod_{k=1}^n \det \zeta_{2k-1} \notag\\ 
&= (\det U_1) (\det V_2)\cdots (\det U_{2n-2}) (\det V_{2n-2}) (\det U_{2n-1}) ,
\end{align}
which 
gives the second identity in \eqref{finalchange}.
\hfill $\Box $

\subsubsection{Second proof of  Lemma \ref{crucial} 
via  Szeg{\H o}'s mapping}
It was tempting to extend to the 
matrix case the method used in the scalar 
one for the Jacobi ensemble. 
The main steps use successively:
\begin {itemize}
\item 
the inverse Szeg{\H o} mapping to turn the problem  on $[0,1]$ into a problem on the unit circle,  
\item  the correspondence between orthogonal polynomials on the unit circle and on the real line,
\item the  Szeg{\H o} recursion for polynomials on the unit circle.
\end{itemize}
To begin with, we transfer the measure on $[0,1]$ to a measure on $[-2, 2]$ by the mapping $x \mapsto 2 -4x$. The new Jacobi matrix $\hat J_n$ is deduced from the original matrix $J_n$ by
\begin{align}
\label{Omega}
\hat J_n =2I_n - 4 \Omega J_n \Omega\,,
\end{align}
where $\Omega$ is a diagonal matrix with alternating blocks $\pm\bf 1$'s on the diagonal.

Let $\hat P_0,\dots ,\hat P_n$ be the monic orthogonal polynomials associated with $\hat J_n$.
From \cite{damanik2008analytic}  Section 2.9 (with reference in particular to \cite{duran1996orthogonal} and \cite{sinap1996orthogonal})
\begin{align}
\det \hat P_n(z) = \det (zI_n- \hat J_n) , 
\end{align}
so that
\begin{align}
\det \hat P_n(z) = \det\left((z-2) I_n + 4 \Omega J_n\Omega\right)= 4^{np} \det\left( \frac{z-2}{4} I_n + J_n\right)
\end{align}
and in particular, 
\begin{align}
\label{detPbold}
\det J_n = 4^{-np} \det \hat P_n (2) , \qquad  \det (I_n -J_n) = (-4)^{-np} \det \hat P_n (-2)\,.
\end{align}
We refer to the definition of the Szeg{\H o} mapping given in Section \ref{susec:Szego}. In this Section, we write $\Sigma_{\mathbb{R}}$ for a matrix measure on the real line and denote by $\Sigma_{\mathbb{T}} = \tilde{\operatorname{Sz}}^{-1}(\Sigma_{\mathbb{R}})$ the preimage under the Szeg{\H o} mapping. 
The correspondence between polynomials orthogonal with respect to $\Sigma_{\mathbb T}$ and with respect to $\Sigma_{\mathbb R}$ 
 is ruled by the following theorem (see Proposition 1 in \cite{Y-M}). It is the 
matrix version of a famous theorem due to Szeg{\H o} \cite{szego1939orthogonal}. Since the notations are slightly different from the usual ones, we rewrite the proof in Section 
\ref{YMproof}. 

\begin{thm}[Yakhlef-Marcell\'an]
\label{Y-M}
Let $\Sigma_{\mathbb R}\in \mathcal{M}_{p,1}([-2,2])$ be a nontrivial matrix measure and denote by $\Sigma_{\mathbb T}=\tilde{\operatorname{Sz}}(\Sigma_{\mathbb{T}})$ the symmetric measure on $\mathbb T$ obtained by the Szeg{\H o} mapping.

If $\hat P_n$  is the $n$-th right monic orthogonal polynomial for $\Sigma_{\mathbb R}$
 and $\Fbold_{2n}$ the $2n$-th right monic orthogonal polynomial\footnote{The right monic OP for $\Sigma_{\mathbb{T}}$ are obtained by applying Gram-Schmidt to $\{\idbold, z\idbold, \dots\}$.} for $\Sigma_{\mathbb T}$, then 
\begin{align}
\hat P_n (z + z^{-1})=   \left[z^{-n} \Fbold_{2n} (z) + z^n \Fbold_{2n}(z^{-1})\right] \tobold^{-1}_n\,,
\end{align}
where
\begin{align}
\label{defto}
\tobold_n := {\bf 1} + \Fbold_{2n}(0) = {\bf 1}- \kbold_{2n-1} \abold_{2n-1} (\kbold_{2n-1})^{-1}
\end{align}
with
\[\kbold_k = \left(\rbold_0\dots\rbold_{k-1}\right)^{-1} , \qquad \rbold_j = ({\bold 1} - \abold_j^2)^{1/2}\,.\]
\end{thm}

From \eqref{detPbold} and \eqref{defto} we deduce 
taking $z= \pm 1$,
\begin{align}
\hat P_n(\pm 2) =  2 
(\pm 1)^n \Fbold_{2n}(\pm 1)\tobold_n^{-1}\,,
\end{align}
hence
\begin{align}
\label{detPbold2}
\det \hat P_n(\pm 2) = 2^p \det (1 - \abold_{2n-1})^{-1} (\pm 1)^{np}\det \Fbold_{2n}(\pm 1)\,.
\end{align}

Recall that the recursion formula  
 expressed for the monic polynomials on the unit circle, in this particular case, is 
\begin{align}
\label{recursionMOPUC2}
z\Fbold_k(z) - \Fbold_{k+1}(z) = z^k \Fbold_k(z^{-1}) \kbold_k \abold_k \kbold_k^{-1}
\end{align}
(see (3.11) in \cite{damanik2008analytic}),
so that 
\[\Fbold_{2n}(1) = \prod_{j=0}^{2n-1} \left({\bf 1} -  \kbold_j \abold_j \kbold_j^{-1}\right)
 , \qquad \Fbold_{2n}(-1) =  \prod_{j=0}^{2n-1} \left({\bf 1} + (-1)^j  \kbold_j \abold_j \kbold_j^{-1}\right)\]
and then
\begin{align}
\det \Fbold_{2n}(1) = \prod_{j=0}^{2n-1} \det \left({\bf 1} -  \abold_j \right)  , \qquad  \det \Fbold_{2n}(-1) =  \prod_{j=0}^{2n-1} \det \left({\bf 1} + (-1)^j  \abold_j \right)\,.
\end{align}

These relations are the matrix extension of Lemma 5.2 of \cite{Killip1}.
Coming back to \eqref{detPbold} and \eqref{detPbold2}, we get
\begin{align}
\det (J_n ) =  2^{-(2n-1)p} \prod_{j=0}^{2n-2}\det({\bf 1} -\abold_j) , \qquad \det (I_n -J_n) = 2^{-(2n-1)p} \prod_{j=0}^{2n-2} \det({\bf 1} + (-1)^j \abold_j) .
\end{align}
The connection with the canonical moments follows then from \eqref{DeWag}. Note that this identity 
still holds  if $\Sigma$ is not nontrivial, as long as ${\bf 0} <\mathcal{U}_k <{\bf 1}$, or equivalently $-{\bf 1} <\abold_{k-1} <{\bf 1}$.

\subsubsection{Proof of Theorem \ref{Y-M}}
\label{YMproof}

In the scalar case, the proof  is given in \cite{simon1} Theorem 13.1.5 or in \cite{simon2} Theorem 1.9.1, with references therein. In the matrix case, one can follow the same scheme.

Since $\Sigma_{\mathbb T}$ is invariant, the Verblunsky coefficients are Hermitian (see Lemma 4.1 in \cite{damanik2008analytic}). 
The matrix Laurent polynomial 
$z^{-n}\Fbold_{2n}(z)+ z^n \Fbold_{2n}(z^{-1})$
 is invariant by $z \mapsto z^{-1}$. Hence there exists a matrix polynomial $\tilde Q_n$ of degree $n$, such that
\begin{align}\label{zbarz} 
z^{-n}\Fbold_{2n}(z)+ z^n \Fbold_{2n}(z^{-1})
 = \tilde Q_n (z + z^{-1})\,,\end{align}
(see for instance Lemma 13.4.2 in \cite{simon2}). 
Collecting terms with highest degrees, we 
 have
\[\tilde Q_n (z+z^{-1}) =
\left(z^n + z^{-n}\right) \tobold_n + \cdots\]
and then
\begin{align}
\label{defQn}
\tilde Q_n (z + z^{-1}) \tobold_n^{-1}  = Q_n (z + z^{-1})
\end{align}
where now $Q_n(x)$ is a monic polynomial of degree $n$. 
Now, let us check that the $\tilde Q_k$ (hence $Q_k$) are orthogonal polynomials for $\Sigma_{\mathbb{R}}$. 
First notice that
\[\tilde Q_k(z+z^{-1}) = z^{-k} \left(\Fbold_{2k}(z) + z^{2k} \Fbold_{2k}( z^{-1})\right)\,.\] 
From 
the Szeg{\H o} mapping and (\ref{defQn}), orthogonality of $\tilde Q_n$ and $\tilde Q_r$ (for $n\not= r$) with respect to $\Sigma_{\mathbb{R}}$ is equivalent to orthogonality (with respect to $\Sigma_{\mathbb T}$) of  
$\Fbold_{2n}(z) + z^{2n} \Fbold_{2n}(z^{-1})$
 and 
$H$ where
\[H(z) = z^{n-r}\left[\Fbold_{2r}(z) + z^{2r} \Fbold_{2r}( z^{-1})\right]\,,\]
which is a polynomial of degree $n+r$ without constant term.
By definition, $\Fbold_{2n}$ is orthogonal to $ z^j{\bf 1}$ for all $j=0, \dots, 2n-1$. Besides, 
$z^{2n} \Fbold_{2n}(z^{-1})$ is (right) orthogonal to
$z^j{\bf 1}$ for $j =1, \dots, 2n$. Indeed, 
\begin{align*}
\int \left[z^{2n} \Fbold_{2n}(\bar z)\right]^\dagger d\Sigma_{\mathbb{T}}(z) z^j &= \int \Fbold_{2n}(\bar z)^\dagger  d\Sigma_{\mathbb{T}}(z) z^{j-2n}\\
&= \int \Fbold_{2n}(z)^\dagger d\Sigma_{\mathbb{T}}(z) z^{2n-j}
\end{align*}
(by invariance of $\Sigma_{\mathbb{T}}$)
and this last integral is ${ \bf 0}$ for $1 \leq j \leq 2n$ due to the orthogonality of $\Fbold_{2n}$ with polynomials of degree at most $2n-1$.

One can then conclude that
$\Fbold_{2n}(z) + z^{2n} \Fbold_{2n}(z^{-1})$
 is orthogonal to $z^k$ for $1\leq k\leq 2n-1$, 
and so to $H$. Summarizing,  
 the $Q_n$'s are the monic polynomials orthogonal with respect to $\Sigma_{\mathbb{R}}$, and then 
$\hat P_n = Q_n$ for every $n$, 
or in other words, by (\ref{zbarz}) and (\ref{defQn})
\begin{align}
\hat P_n (z+z^{-1}) =  
\left[z^{-n}\Fbold_{2n}(z) + z^{n} \Fbold_{2n}(z^{-1})\right]\tobold_n^{-1}\,.
\end{align}
\hfill $ \Box $

\bibliographystyle{plain}
\bibliography{bibclean}

\end{document}